\documentclass[11pt,a4paper]{article}
\date{}
\makeatletter
\@addtoreset{equation}{section}
\makeatother

\usepackage[colorlinks,
            linkcolor=blue,
            anchorcolor=blue,
            citecolor=blue]{hyperref}
\usepackage{hyperref}
\hypersetup{colorlinks=true,
            linkcolor=blue,
            anchorcolor=blue,
            citecolor=blue}
\usepackage{hyperref,amsmath,amssymb,amscd}
\usepackage{amssymb,amsmath,enumerate}
\usepackage{indentfirst}
\usepackage{url}
\usepackage[capitalise]{cleveref}
\usepackage{amsthm}
\usepackage{amsmath}
\AtBeginDocument{\hypersetup{pdfborder={0 0 0.01}}}
\newtheorem{theorema}{Theorem}[section]

\newtheorem{thm}{Theorem}[section]
\newtheorem{cor}[thm]{Corollary}
\newtheorem{lem}[thm]{Lemma}

\theoremstyle{definition}

\newtheorem{rem}[thm]{\bf Remark}
\allowdisplaybreaks[4]

\usepackage{enumitem}
\usepackage{extarrows}
\usepackage{amsthm}
\usepackage{setspace}
\hypersetup{urlcolor=blue, citecolor=red}
\usepackage{cite}

\begin{document}

\title{\bf
Stability estimates for $L^p$-Caffarelli-Kohn-Nirenberg inequalities\footnote{Supported by National Natural Science Foundation of China (No. 12371120).}}
\author{{Xiao-Ping Chen,\ \ Chun-Lei Tang\footnote{Corresponding author.\newline
\indent\,\,\, \emph{E-mail address:} xpchen\_maths@163.com (X.-P. Chen); tangcl@swu.edu.cn (C.-L. Tang).}}\\
{\small \emph{School of Mathematics and Statistics, Southwest University,  Chongqing {\rm400715},}}\\
{\small \emph{People's Republic of China}}}
\maketitle
\baselineskip 17pt

\noindent {\bf Abstract}:\ Based on some new vector inequalities established by Figalli and Zhang [\emph{Duke Math. J.} \textbf{171} (2022),  2407--2459], we study the stability of the scale invariant and the scale non-invariant $L^p$-Caffarelli-Kohn-Nirenberg inequalities, which fills the recent work of Do \emph{et al.} [$L^p$-Caffarelli-Kohn-Nirenberg inequalities and their stabilities, arXiv: 2310.07083] for $1<p<2$, and also extends some results of Cazacu \emph{et al.} [\emph{J. Math. Pures Appl. (9)} \textbf{182} (2024), 253--284] to a general case for $L^p$-Caffarelli-Kohn-Nirenberg inequalities with $1<p<N$.

\vspace{0.25em}


\vspace{0.25em}


\section{Introduction and main results}

\noindent In this paper, we focus on the stability of several classes of $L^p$-Caffarelli-Kohn-Nirenberg ($L^p$-CKN for short) inequalities. The novelty of this paper is stated as follows.
\begin{enumerate}
[itemsep=0pt, topsep=2pt, parsep=0pt]

\item[$(1)$]
Here we emphasize that the stability of the scale invariant and the scale non-invariant $L^p$-CKN inequalities for $1<p<2$ is novel. Furthermore, our main results improve and generalize those of \cite{Cazacu24,Do23}.

\item[$(2)$]
A nonexistence result about the stability of the scale invariant $L^p$-CKN inequalities has also been discussed.

\item[$(3)$]
As a byproduct, we prove a series of weighted  $L^p$-Poincar\'{e} inequalities (including the case $p=2$) for the log-concave measure on the Borel sets of $\mathbb{R}^N$.

\end{enumerate}

We provide some related results in the first part of this section about sharp constants, extremizers and the stability of $L^p$-CKN inequalities (including $p=2$). Our main results are described in the second part of this section. We present the main novelty, difficulties and strategy of this paper in the third part of this section. The final part of this section contains the outline for the remainder of this paper.

\subsection{Overview and motivation}\label{sub-1.1}

\noindent In 1984, Caffarelli, Kohn and Nirenberg  \cite{Caffarelli84} first introduced Caffarelli-Kohn-Nirenberg (CKN for short) inequalities. Afterwards, CKN inequalities (which include Sobolev inequalities, Hardy-Sobolev inequalities, Gagliardo-Nirenberg inequalities and so on) had been studied extensively because of their significance in various aspects of mathematics, such as \cite{Catrina09,Catrina01,Cazacu21,Cazacu23,
Cazacu24,Costa08,Dong18,Lam17,Mallick22,Nguyen15} and the reference therein.

Especially, the following $L^2$-CKN inequalities are important subfamilies of CKN inequalities,
\begin{align}\label{2-ckn-ineq}
&\left(\int_{\mathbb{R}^N}
\frac{|\nabla u|^2}
{|x|^{2b}}\mathrm{d}x
\right)^{\frac{1}{2}}
\left(\int_{\mathbb{R}^N}
\frac{|u|^2}{|x|^{2a}}
\mathrm{d}x\right)^{\frac{1}{2}}
\geq C(N,a,b)
\int_{\mathbb{R}^N}
\frac{|u|^2}{|x|^{a+b+1}}
\mathrm{d}x,
\end{align}
for all $u\in\mathcal{C}_0^\infty(\mathbb{R}^N
\setminus\{\mathbf{0}\})$; and $\mathcal{C}_0^\infty
(\mathbb{R}^N\setminus\{\mathbf{0}\})$ consists of smooth functions compactly supported on $\mathbb{R}^N\setminus\{\mathbf{0}\}$. Motivated by \cite[Theorem 4.4]{Do23}, we present a more general $L^p$-CKN inequalities,
\begin{align}\label{2-ckn-ineq-p}
&\left(\int_{\mathbb{R}^N}
\frac{|\nabla u|^p}
{|x|^{pb}}\mathrm{d}x
\right)^{\frac{1}{p}}
\left(\int_{\mathbb{R}^N}
\frac{|u|^p}{|x|^{pa}}
\mathrm{d}x\right)^{\frac{p-1}{p}}
\geq C(N,p,a,b)
\int_{\mathbb{R}^N}
\frac{|u|^p}{|x|^{(p-1)a+b+1}}
\mathrm{d}x.
\end{align}

\subsubsection{Sharp constants and extremizers for $L^p$-CKN inequalities (including $p=2$)}

\noindent On the one hand, we describe the research status of the $L^2$-CKN inequality \eqref{2-ckn-ineq} about the sharp constants and extremizers. Costa \cite{Costa08} first analyzed the constant $C(N,a,b)$ of \eqref{2-ckn-ineq}, the author obtained that $C(N,a,b)=\frac{|N-(a+b+1)|}{2}$ is the sharp constant of \eqref{2-ckn-ineq}. Later, Catrina and Costa \cite{Catrina09} defined the following regions to study the sharp constants and extremizers of \eqref{2-ckn-ineq},
\begin{equation}\label{1.3}
\begin{cases}
\mathcal{P}_1:=\left\{(a,b):b-a+1>0,\ b\le\frac{N-2}{2}\right\};\\[0.5mm]
\mathcal{P}_2:=\left\{(a,b):b-a+1<0,\ b\ge\frac{N-2}{2}\right\};\\[0.5mm]
\mathcal{P}:=\mathcal{P}_1\cup\mathcal{P}_2;\\[0.5mm]
\mathcal{Q}_1:=\left\{(a,b):b-a+1<0,\ b\le\frac{N-2}{2}\right\};\\[0.5mm]
\mathcal{Q}_2:=\left\{(a,b):b-a+1>0,\ b\ge\frac{N-2}{2}\right\};\\[0.5mm]
\mathcal{Q}:=\mathcal{Q}_1\cup\mathcal{Q}_2.
\end{cases}
\end{equation}
The main results of \cite{Catrina09} state as follows.

\begin{theorema}
[{\!\rm{\!\cite[Theorem 1]{Catrina09}}}]
\label{thm-a}
For a nonzero constant $\alpha$, there hold the following results.
\begin{enumerate}
[itemsep=0pt, topsep=2pt, parsep=0pt]

\item[(1)]
If $(a,b)\in\mathcal{P}$, the sharp constant $C(N,a,b)=\frac{|N-(a+b+1)|}{2}$ is achieved by the functions $u(x)=\alpha e^
    {\frac{\beta|x|^{b-a+1}}{b-a+1}},$
    where $\beta<0$ in $\mathcal{P}_1$ and $\beta>0$ in $\mathcal{P}_2$.

\item[(2)]
If $(a,b)\in\mathcal{Q}$, the sharp constant $C(N,a,b)=\frac{|N-(3b-a+3)|}{2}$ is achieved by the functions $u(x)=\alpha|x|^{2(b+1)-N}
    e^{\frac{\beta|x|^{b-a+1}}{b-a+1}},$
    where $\beta>0$ in $\mathcal{Q}_1$ and $\beta<0$ in $\mathcal{Q}_2$.

\item[(3)]
If $a=b+1$, the sharp constant $C(N,b+1,b)=\frac{|N-2(b+1)|}{2}$ cannot be achieved.
\end{enumerate}
\end{theorema}

\noindent Recently, Cazacu, Flynn and Lam in   \cite{Cazacu21} provided a short proof of the following  inequalities, for all $u\in\mathcal{C}_0^\infty(\mathbb{R}^N
\setminus\{\mathbf{0}\})$,
\begin{align}\label{2-ckn-re-ineq}
&\left(\int_{\mathbb{R}^N}
\frac{|x\cdot\nabla u|^2}
{|x|^{2b+2}}\mathrm{d}x
\right)^{\frac{1}{2}}
\left(\int_{\mathbb{R}^N}
\frac{|u|^2}{|x|^{2a}}
\mathrm{d}x\right)^{\frac{1}{2}}
\geq \tilde{C}(N,a,b)
\int_{\mathbb{R}^N}
\frac{|u|^2}{|x|^{a+b+1}}
\mathrm{d}x,
\end{align}
where $\tilde{C}(N,a,b)=C(N,a,b)>0$ is the sharp constant of \eqref{2-ckn-re-ineq}. Moreover, instead of requiring $\nabla u$, \eqref{2-ckn-re-ineq} only needs the radial derivative $\partial_r u=\frac{x}{|x|}\cdot\nabla u$.

On the other hand, we consider the sharp constants and extremizers of the $L^p$-CKN inequalities \eqref{2-ckn-ineq-p}. In the recent work \cite{Do23}, Do, Flynn, Lam and Lu obtained the following results.

\begin{theorema}
[{\!\rm{\!\cite[Corollary 1.2]{Do23}}}]
\label{thm-b}
Let $N\ge1$ and $p>1$, for all $u\in\mathcal{C}_0^\infty
(\mathbb{R}^N\setminus\{\mathbf{0}\})$, there holds
\begin{equation}\label{ckn}
\left(\int_{\mathbb{R}^N}
\frac{|\nabla u|^p}
{|x|^{pb}}\mathrm{d}x
\right)^{\frac{1}{p}}
\left(\int_{\mathbb{R}^N}
\frac{|u|^p}{|x|^{pa}}
\mathrm{d}x\right)^{\frac{p-1}{p}}
\geq \frac{\left|N-(p-1)a-b-1\right|}{p}
\int_{\mathbb{R}^N}
\frac{|u|^p}{|x|^{(p-1)a+b+1}}
\mathrm{d}x.
\end{equation}
\noindent Moreover,
\begin{enumerate}
[itemsep=0pt, topsep=2pt, parsep=0pt]

\item[(1)]
if $b-a+1>0$ and $b\le\frac{N-p}{p}$, the sharp constant $\frac{N-(p-1)a-b-1}{p}$ is only attained by the functions $u(x)=\alpha e^{\frac{\beta}{b-a+1}|x|^{b-a+1}}$ with $\alpha\in\mathbb{R}$ and $\beta<0$;

\item[(2)]
if $b-a+1<0$ and $b\ge\frac{N-p}{p}$, the sharp constant $\frac{(p-1)a+b+1-N}{p}$ is only attained by the functions $u(x)=\alpha e^{\frac{\beta}{b-a+1}|x|^{b-a+1}}$ with $\alpha\in\mathbb{R}$ and $\beta>0$.

\end{enumerate}
\end{theorema}

For more results about sharp constants and extremizers of second order CKN inequalities, we refer the interested readers to \cite{Cazacu22,Cazacu23,Chen24} and the reference therein.

\subsubsection{The stability of $L^p$-CKN inequalities (including $p=2$)}\label{sub-1.1.2}

\noindent The stability of geometric inequalities got a growing attention after a question proposed by Brezis and Lieb in \cite{Brezis85}. Subsequently, Bianchi and Egnell \cite{Bianchi91} gave an affirmative answer to this question. For more stability results of geometric inequalities, see for example \cite{Bartsch03,Carlen13,Cazacu24,Chen23,
Chen13,Chen25,Cianchi09,De23,Do23,Figalli13,
Lu00,Nguyen19} and the reference therein.

Now, we recall some stability results of $L^2$-CKN inequalities. Recently, in \cite{Cazacu24}, Cazacu, Flynn, Lam and Lu considered the following scale invariant $L^2$-CKN inequalities:
\begin{equation}\label{2-ckn}
\left(\int_{\mathbb{R}^N}
\frac{|\nabla u|^2}
{|x|^{2b}}\mathrm{d}x
\right)^{\frac{1}{2}}
\left(\int_{\mathbb{R}^N}
\frac{|u|^2}{|x|^{2a}}
\mathrm{d}x\right)^{\frac{1}{2}}
\ge\left|\frac{N-(a+b+1)}{2}\right|
\int_{\mathbb{R}^N}
\frac{|u|^2}{|x|^{a+b+1}}
\mathrm{d}x.
\end{equation}
If $0\le b\le\frac{N-2}{2}$, $a\le\frac{Nb}{N-2}$ and $a+b+1=\frac{2bN}{N-2}$ (which belongs to $\mathcal{P}_1$), they obtained that
\begin{align*}
&\left(\int_{\mathbb{R}^N}
\frac{|\nabla u|^2}
{|x|^{2b}}\mathrm{d}x
\right)^{\frac{1}{2}}
\left(\int_{\mathbb{R}^N}
\frac{|u|^2}{|x|^{2a}}
\mathrm{d}x\right)^{\frac{1}{2}}
-\frac{N-(a+b+1)}{2}
\int_{\mathbb{R}^N}
\frac{|u|^2}{|x|^{a+b+1}}
\mathrm{d}x
\nonumber\\&\quad
\ge C(N,a,b)\inf_{c\in\mathbb{R},\lambda>0}
\int_{\mathbb{R}^N}
\frac{\left|u
-ce^
{-\frac{\lambda|x|^{b-a+1}}{b-a+1}}\right|^2}
{|x|^{a+b+1}}
\mathrm{d}x.
\end{align*}
Moreover, they also considered the stability of the following $L^2$-CKN inequalities:
\begin{align}\label{1.7}
&\int_{\mathbb{R}^N}
\frac{|\nabla u|^2}{|x|^{2b}}\mathrm{d}x
+\int_{\mathbb{R}^N}
\frac{|u|^2}{|x|^{2a}}\mathrm{d}x
\ge\left|a+b+1-N\right|
\int_{\mathbb{R}^N}
\frac{|u|^2}{|x|^{a+b+1}}\mathrm{d}x.
\end{align}
If $0\le b<\frac{N-2}{2}$ and $a\le\frac{Nb}{N-2}$, there holds
\begin{align*}
&\int_{\mathbb{R}^N}
\frac{|\nabla u|^2}{|x|^{2b}}\mathrm{d}x
+\int_{\mathbb{R}^N}
\frac{|u|^2}{|x|^{2a}}\mathrm{d}x
-\left[N-\left(a+b+1\right)\right]
\int_{\mathbb{R}^N}
\frac{|u|^2}{|x|^{a+b+1}}\mathrm{d}x
\\&\quad
\ge C(N,a,b)\inf_{c\in\mathbb{R}}
\int_{\mathbb{R}^N}
\frac{\left|u
-ce^
{-\frac{|x|^{b-a+1}}{b-a+1}}\right|^2}
{|x|^{\frac{2bN}{N-2}}}
\mathrm{d}x.
\end{align*}
More importantly, if $0\le b<\frac{N-2}{2}$, $a\le\frac{Nb}{N-2}$ and $a+b+1=\frac{2bN}{N-2}$, a stronger stability version of \eqref{1.7} had also been obtained in \cite{Cazacu24}, for some $C=C(N,a,b)>0$,
\begin{align*}
&\int_{\mathbb{R}^N}
\frac{|\nabla u|^2}{|x|^{2b}}\mathrm{d}x
+\int_{\mathbb{R}^N}
\frac{|u|^2}{|x|^{2a}}\mathrm{d}x
-\left[N-\left(a+b+1\right)\right]
\int_{\mathbb{R}^N}
\frac{|u|^2}{|x|^{a+b+1}}\mathrm{d}x
\\&\quad
\ge C\inf_{c\in\mathbb{R}}
\left\{
\int_{\mathbb{R}^N}\!
\frac{\left|\nabla\!\left(\!u
\!-\!ce^{-\frac{|x|^{b-a+1}}
{b-a+1}}\!\right)\!\right|^2}
{|x|^{2b}}
\mathrm{d}x
\!+\!\int_{\mathbb{R}^N}\!
\frac{\left|u
\!-\!ce^{-\frac{|x|^{b-a+1}}
{b-a+1}}\right|^2}
{|x|^{2a}}
\mathrm{d}x
\!+\!\int_{\mathbb{R}^N}\!
\frac{\left|u
\!-\!ce^
{-\frac{|x|^{b-a+1}}{b-a+1}}\right|^2}
{|x|^{a+b+1}}\mathrm{d}x\right\}.
\end{align*}

After this, the authors in \cite{Do23} studied the stability of the following $L^2$-CKN inequalities:
\begin{equation}\label{2-ckn-2}
\left(\int_{\mathbb{R}^N}
\frac{|\nabla u|^2}
{|x|^{2b}}\mathrm{d}x
\right)^{\frac{1}{2}}
\left(\int_{\mathbb{R}^N}
\frac{|u|^2}{|x|^{2a}}
\mathrm{d}x\right)^{\frac{1}{2}}
\ge\left|\frac{N-(3b-a+3)}{2}\right|
\int_{\mathbb{R}^N}
\frac{|u|^2}{|x|^{a+b+1}}
\mathrm{d}x.
\end{equation}
If $\frac{N-2}{2}<b\le N-2$ and $N(b-a+3)=2(3b-a+3)$ (which belongs to $\mathcal{Q}_2$), they proved that
\begin{align*}
&\left(\int_{\mathbb{R}^N}
\frac{|\nabla u|^2}
{|x|^{2b}}\mathrm{d}x
\right)^{\frac{1}{2}}
\left(\int_{\mathbb{R}^N}
\frac{|u|^2}{|x|^{2a}}
\mathrm{d}x\right)^{\frac{1}{2}}
-\frac{(3b-a+3)-N}{2}
\int_{\mathbb{R}^N}
\frac{|u|^2}{|x|^{a+b+1}}
\mathrm{d}x
\nonumber\\&\quad
\ge C(N,a,b)
\inf_{c\in\mathbb{R},\lambda>0}
\int_{\mathbb{R}^N}
\frac{\left|u
-c|x|^{2b+2-N}e^
{-\frac{\lambda|x|^{b-a+1}}{b-a+1}}\right|^2}
{|x|^{a+b+1}}
\mathrm{d}x.
\end{align*}

Our main goal of this paper is to study the stability of $L^p$-CKN inequalities. To our knowledge, there are very few results about it, see  Do \emph{et al.} \cite{Do23}. In fact, they proved the following result.

\begin{theorema}
[{\!\rm{\!\cite[Theorem 1.9]{Do23}}}]
\label{thm-c}
Let $p\ge2$, $0\le b<\frac{N-p}{p}$, $a<\frac{Nb}{N-p}$ and $(p-1)a+b+1=\frac{pbN}{N-p}$, for each $u\in\mathcal{C}_0^\infty
(\mathbb{R}^N\setminus\{\mathbf{0}\})$, there is a universal constant $C(N,p,b)>0$ such that
\begin{align*}
&\left(\int_{\mathbb{R}^N}
\frac{|\nabla u|^p}
{|x|^{pb}}\mathrm{d}x
\right)^{\frac{1}{p}}
\left(\int_{\mathbb{R}^N}
\frac{|u|^p}{|x|^{pa}}
\mathrm{d}x\right)^{\frac{p-1}{p}}
-\frac{N-(p-1)a-b-1}{p}
\int_{\mathbb{R}^N}
\frac{|u|^p}{|x|^{(p-1)a+b+1}}
\mathrm{d}x
\\&\quad\ge
C(N,p,b)
\inf_{c\in\mathbb{R},\lambda>0}
\int_{\mathbb{R}^N}
\frac{\left|u
-ce^
{-\frac{\lambda|x|^{b-a+1}}{b-a+1}}\right|^p}
{|x|^{(p-1)a+b+1}}
\mathrm{d}x.
\end{align*}
\end{theorema}

Motivated by the results shown above, it seems reasonable to ask some questions.

\vspace{0.5em}

\textbf{$\bullet$\,\emph{Question 1}}: \emph{does there exist  a stability version of \eqref{ckn} when $1<p<2$?}

\textbf{$\bullet$\,\emph{Question 2:}} \emph{when $N\ge1$, $p>1$, we wonder if the stability of the following scale non-invariant $L^p$-CKN inequalities exists,}
\begin{align}\label{sniLpI}
\int_{\mathbb{R}^N} \frac{|\nabla u|^p}{|x|^{pb}}\mathrm{d}x +(p-1)\int_{\mathbb{R}^N} \frac{|u|^p}{|x|^{pa}}\mathrm{d}x \ge\left|N-(p-1)a-b-1\right|\int_{\mathbb{R}^N} \frac{|u|^p}{|x|^{(p-1)a+b+1}} \mathrm{d}x.
\end{align}

\vspace{0.5em}

\noindent We give an affirmative answer to Question 1 (resp. Question 2) in Theorem \ref{thm-1} (resp. Theorems \ref{thm-3}, \ref{thm-4} and \ref{thm-5}).

\subsection{Main results}

\noindent Now, let us present our main results.

\subsubsection{The stability results for the scale invariant $L^p$-CKN inequalities}

\noindent In this subsection, we provide the stability results for the scale invariant $L^p$-CKN inequalities. Let
\[
\mathcal{M}_{a,b}:
=\left\{\alpha
e^{\frac{\beta}{b-a+1}|x|^{b-a+1}}:
\alpha\in\mathbb{R},\
\frac{\beta}{b-a+1}<0\right\}
\]
be the set of extremal functions of \eqref{ckn}.

For the Sobolev inequality, Figalli and Zhang in \cite[Theorem 1.1]{Figalli22} proved quantitative stability results in terms of  $\|\nabla (u-v^*)\|_{L^p(\mathbb{R}^N)}$ for a function $u$ and its corresponding closest extremal function $v^*$ of the Sobolev inequality. We also expect similar results hold for the scale invariant $L^p$-CKN inequality \eqref{ckn}. However, motivated by \cite[Proposition 1.7]{McCurdy21}, we prove that there exists no quantitative stability results in terms of gradient for the scale invariant $L^p$-CKN inequality \eqref{ckn} in the following theorem. This indicates that Theorem \ref{thm-c} is sharp.

\begin{thm}\label{thm-6}
Let $N\ge1$, $p>1$, $b-a+1>0$ and $b\le\frac{N-p}{p}$. Then for any two nonnegative constants $C_1$ and $C_2$ such that $C_1+C_2>0$, there exists a function $u\in\mathcal{C}_0^\infty
(\mathbb{R}^N\setminus
\{\mathbf{0}\})\setminus\mathcal{M}_{a,b}$ such that the following estimate holds,
\begin{align*}
&\left(\int_{\mathbb{R}^N}
\frac{|\nabla u|^p}{|x|^{pb}}\mathrm{d}x
\right)^{\frac{1}{p}}
\left(\int_{\mathbb{R}^N}
\frac{|u|^p}{|x|^{pa}}
\mathrm{d}x\right)^{\frac{p-1}{p}}
-\frac{N-(p-1)a-b-1}{p}
\int_{\mathbb{R}^N}
\frac{|u|^p}{|x|^{(p-1)a+b+1}}
\mathrm{d}x
\\&\quad\le
C_1\inf_{v\in\mathcal{M}_{a,b}}\int_{\mathbb{R}^N}
\frac{|\nabla (u-v)|^p}{|x|^{pb}}\mathrm{d}x
+
C_2
\inf_{v\in\mathcal{M}_{a,b}}
\int_{\mathbb{R}^N}
\frac{
\left|u-v\right|^p}{|x|^{pa}}
\mathrm{d}x.
\end{align*}
\end{thm}

After completing this paper, it was found that there exist a similar result to Theorem \ref{thm-6} in \cite[Proposition 1.7]{Zhang24}. For the completeness of the paper, we retain this result here.

Theorem \ref{thm-6} is very useful and effective for our further purposes. In view of Theorem \ref{thm-6}, we provide the stability results of the scale invariant $L^p$-CKN inequality \eqref{ckn} in the following two theorems, precisely for the case $1<p<2$ and for the case $p\ge2$.

\begin{thm}\label{thm-1}
Let $N>2$, $1<p<2$, $0\le b
<\frac{N-p}{p}$ and
$a=\frac{Nb}{N-p}$. For each $u\in\mathcal{C}_0^\infty
(\mathbb{R}^N\setminus\{\mathbf{0}\})\setminus\{0\}$,
there is a universal constant $C(N,p,b)>0$ such that
\begin{align*}
&\left(\int_{\mathbb{R}^N}
\frac{|\nabla u|^p}{|x|^{pb}}\mathrm{d}x
\right)^{\frac{1}{p}}
\left(\int_{\mathbb{R}^N}
\frac{|u|^p}{|x|^{pa}}
\mathrm{d}x\right)^{\frac{p-1}{p}}
-\frac{N-(p-1)a-b-1}{p}
\int_{\mathbb{R}^N}
\frac{|u|^p}{|x|^{(p-1)a+b+1}}
\mathrm{d}x
\\&\quad\ge
C(N,p,b)\left(
\frac{\int_{\mathbb{R}^N}
\frac{|\nabla u|^p}{|x|^{pb}}\mathrm{d}x}
{\int_{\mathbb{R}^N}
\frac{|u|^p}{|x|^{pa}}\mathrm{d}x}
\right)^{\frac{1}{p}}
\inf_{c\in\mathbb{R},\lambda>0}
\int_{\mathbb{R}^N}
\frac{
\left||u|^{\frac{p-2}{2}}u
-|c|^{\frac{p-2}{2}}c
e^{-\frac{\frac{p}{2}|x|^{b-a+1}}
{(b-a+1)\lambda^{b-a+1}}}
\right|^2}{|x|^{pa}}
\mathrm{d}x,
\end{align*}
equivalently,
\begin{align*}
1-\frac{\frac{N-(p-1)a-b-1}{p}
\int_{\mathbb{R}^N}
\frac{|u|^p}{|x|^{(p-1)a+b+1}}
\mathrm{d}x}{
\left({\int_{\mathbb{R}^N}
\frac{|\nabla u|^p}{|x|^{pb}}\mathrm{d}x}
\right)^{\frac{1}{p}}
\left(\int_{\mathbb{R}^N}
\frac{|u|^p}{|x|^{pa}}\mathrm{d}x
\right)^{\frac{p-1}{p}}}
\ge
C(N,p,b)
\inf_{v\in\mathcal{M}_{a,b}}
\left(\frac{\int_{\mathbb{R}^N}\frac{
\left||u|^{\frac{p-2}{2}}u
-|v|^{\frac{p-2}{2}}v\right|^2}{|x|^{pa}}
\mathrm{d}x}{\int_{\mathbb{R}^N}
\frac{|u|^p}{|x|^{pa}}\mathrm{d}x}\right).
\end{align*}
\end{thm}

\begin{thm}\label{thm-2}
Let $N>p\ge2$, $0\le b
<\frac{N-p}{p}$ and
$a=\frac{Nb}{N-p}$.
For all $u\in\mathcal{C}_0^\infty
(\mathbb{R}^N\setminus
\{\mathbf{0}\})\setminus\{0\}$, there is a universal constant $C(N,p,b)>0$ such that
\begin{align*}
&\left(\int_{\mathbb{R}^N}
\frac{|\nabla u|^p}{|x|^{pb}}\mathrm{d}x
\right)^{\frac{1}{p}}
\left(\int_{\mathbb{R}^N}
\frac{|u|^p}{|x|^{pa}}
\mathrm{d}x\right)^{\frac{p-1}{p}}
-\frac{N-(p-1)a-b-1}{p}
\int_{\mathbb{R}^N}
\frac{|u|^p}{|x|^{(p-1)a+b+1}}
\mathrm{d}x
\\&\quad\ge
C(N,p,b)\left(
\frac{\int_{\mathbb{R}^N}
\frac{|\nabla u|^p}{|x|^{pb}}\mathrm{d}x}
{\int_{\mathbb{R}^N}
\frac{|u|^p}{|x|^{pa}}\mathrm{d}x}
\right)^{\frac{1}{p}}
\inf_{c\in\mathbb{R},\lambda>0}
\int_{\mathbb{R}^N}
\frac{
\left|u
-c
e^{-\frac{|x|^{b-a+1}}
{(b-a+1)\lambda^{b-a+1}}}
\right|^p}{|x|^{pa}}
\mathrm{d}x,
\end{align*}
which is equivalent to
\begin{align*}
1-\frac{\frac{N-(p-1)a-b-1}{p}
\int_{\mathbb{R}^N}
\frac{|u|^p}{|x|^{(p-1)a+b+1}}
\mathrm{d}x}{
\left({\int_{\mathbb{R}^N}
\frac{|\nabla u|^p}{|x|^{pb}}\mathrm{d}x}
\right)^{\frac{1}{p}}
\left(\int_{\mathbb{R}^N}
\frac{|u|^p}{|x|^{pa}}\mathrm{d}x
\right)^{\frac{p-1}{p}}}
\ge C(N,p,b)
\inf_{v\in\mathcal{M}_{a,b}}
\left(\frac{\int_{\mathbb{R}^N}
\frac{
\left|u-v\right|^p}{|x|^{pa}}
\mathrm{d}x}{\int_{\mathbb{R}^N}
\frac{|u|^p}{|x|^{pa}}\mathrm{d}x}\right).
\end{align*}
\end{thm}

\begin{rem}
Theorem \ref{thm-6} indicates that there is no quantitative stability version of the scale invariant $L^p$-CKN inequality \eqref{ckn} when adding $\inf_{v\in\mathcal{M}_{a,b}}\int_{\mathbb{R}^N}
\frac{|\nabla (u-v)|^p}{|x|^{pb}}\mathrm{d}x$ or
$\inf_{v\in\mathcal{M}_{a,b}}\int_{\mathbb{R}^N} \frac{|u-v|^p}{|x|^{pa}}\mathrm{d}x$
on the right hand side of \eqref{ckn} for $u\in\mathcal{C}_0^\infty
(\mathbb{R}^N\setminus
\{\mathbf{0}\})\setminus\mathcal{M}_{a,b}$. Thereby, in order to consider the stability result involving $\inf_{v\in\mathcal{M}_{a,b}}\int_{\mathbb{R}^N} \frac{|u-v|^p}{|x|^{pa}}\mathrm{d}x$ of the scale invariant $L^p$-CKN inequality \eqref{ckn}, the term $\left(
\frac{\int_{\mathbb{R}^N}
{|\nabla u|^p}/{|x|^{pb}}\mathrm{d}x}
{\int_{\mathbb{R}^N}
{|u|^p}/{|x|^{pa}}\mathrm{d}x}
\right)^{{1}/{p}}$ need.

Furthermore, compared with Theorem \ref{thm-6}, it is easy to verify that the results of Theorems \ref{thm-1} and \ref{thm-2} are scale invariant.
\end{rem}

As a consequence of Theorem \ref{thm-2}, if $p=2$ and $a=b=0$, the stability result reduces into that of the Hydrogen Uncertainty Principle (a special case of the $L^2$-CKN inequality \eqref{2-ckn}, which can be used to show the stability of a hydrogenic atom in a magnetic field, see \cite{Frohlich86,Lieb10} for details).

\begin{cor}\label{cor-7}
Assume that $N\ge3$. For all $u\in\mathcal{C}_0^\infty
(\mathbb{R}^N\setminus
\{\mathbf{0}\})\setminus\{0\}$, there exists a universal constant $C(N)>0$ such that
\begin{align*}
&\left(\int_{\mathbb{R}^N}
|\nabla u|^2\mathrm{d}x\right)^{\frac{1}{2}}
\left(\int_{\mathbb{R}^N}|u|^2
\mathrm{d}x\right)^{\frac{1}{2}}
-\frac{N-1}{2}
\int_{\mathbb{R}^N}
\frac{|u|^2}{|x|}
\mathrm{d}x
\\&\quad\ge
C(N)\left(
\frac{\int_{\mathbb{R}^N}
|\nabla u|^2\mathrm{d}x}
{\int_{\mathbb{R}^N}|u|^2\mathrm{d}x}
\right)^{\frac{1}{2}}
\inf_{c\in\mathbb{R},\lambda>0}
\int_{\mathbb{R}^N}
\left|u-ce^{-\frac{|x|}
{\lambda}}\right|^2
\mathrm{d}x,
\end{align*}
which is equivalent to
\begin{align*}
1-\frac{\frac{N-1}{2}
\int_{\mathbb{R}^N}
\frac{|u|^2}{|x|}
\mathrm{d}x}{
\left({\int_{\mathbb{R}^N}
|\nabla u|^2\mathrm{d}x}
\right)^{\frac{1}{2}}
\left(\int_{\mathbb{R}^N}
|u|^2\mathrm{d}x
\right)^{\frac{1}{2}}}
\ge C(N)
\inf_{v\in\mathcal{M}_{0,0}}
\left(\frac{\int_{\mathbb{R}^N}
\left|u-v\right|^2
\mathrm{d}x}{\int_{\mathbb{R}^N}
|u|^2\mathrm{d}x}\right).
\end{align*}
\end{cor}

\begin{rem}
As far as we know, for the stability of $L^p$-CKN inequalities, there are few works  \cite{Cazacu24,Do23}. Compared with \cite{Cazacu24,Do23}, some discussions about Theorems \ref{thm-1}, \ref{thm-2} and Corollary \ref{cor-7} are presented below.
\begin{enumerate}[itemsep=0pt, topsep=0pt, parsep=0pt]

\item[(1)]
The result of Theorem \ref{thm-1} is novel for the case $1<p<2$.

\item[(2)]
In contrast with \cite{Do23}, the authors obtained the stability of \eqref{ckn} when adding
\[
\inf_{v\in\mathcal{M}_{a,b}}\int_{\mathbb{R}^N}
\frac{\left|u-v\right|^p}
{|x|^{(p-1)a+b+1}}
\mathrm{d}x
\]
on the right hand side of \eqref{ckn}. Thereby, this paper can be regarded as the complementary work of \cite{Do23}.

\item[(3)]
In view of Corollary \ref{cor-7}, this paper provides another idea to prove the stability of Hydrogen Uncertainty Principle.

\end{enumerate}
\end{rem}

\subsubsection{The stability results for the scale non-invariant $L^p$-CKN inequalities}

\noindent Another purpose of this paper is to show the stability results of the scale non-invariant $L^p$-CKN inequalities.

\begin{thm}\label{thm-3}
Let $N>2$, $1<p<2$, $0\le b
<\frac{N-p}{p}$ and
$a=\frac{Nb}{N-p}$.
For all $u\in\mathcal{C}_0^\infty
(\mathbb{R}^N\setminus
\{\mathbf{0}\})$, there is a universal constant $C(N,p,b)>0$ such that
\begin{align*}
&\int_{\mathbb{R}^N} \frac{|\nabla u|^p}{|x|^{pb}}\mathrm{d}x +(p-1)\int_{\mathbb{R}^N} \frac{|u|^p}{|x|^{pa}}\mathrm{d}x -\left[N-(p-1)a-b-1\right] \int_{\mathbb{R}^N} \frac{|u|^p}{|x|^{(p-1)a+b+1}} \mathrm{d}x
\\&\quad
\ge
C(N,p,b)
\inf_{c\in\mathbb{R}}
\int_{\mathbb{R}^N}
\frac{
\left||u|^{\frac{p-2}{2}}u
-|c|^{\frac{p-2}{2}}c
e^{-\frac{\frac{p}{2}|x|^{b-a+1}}
{b-a+1}}
\right|^2}{|x|^{pa}}
\mathrm{d}x.
\end{align*}
\end{thm}

Compared with Theorem \ref{thm-6}, we now show the difference about stability results between the scale invariant $L^p$-CKN inequality \eqref{ckn} and the scale non-invariant $L^p$-CKN inequality \eqref{sniLpI}. More precisely, there exist stability results  involving gradient terms of $L^p$-CKN inequality \eqref{sniLpI}.

\begin{thm}\label{thm-4}
Let $N>p\ge2$, $0\le b
<\frac{N-p}{p}$ and
$a=\frac{Nb}{N-p}$.
For all $u\in\mathcal{C}_0^\infty
(\mathbb{R}^N\setminus
\{\mathbf{0}\})$, there is a universal constant $C(N,p,b)>0$ such that
\begin{align*}
&\int_{\mathbb{R}^N}
\frac{|\nabla u|^p}{|x|^{pb}}\mathrm{d}x
+(p-1)\int_{\mathbb{R}^N}
\frac{|u|^p}{|x|^{pa}}\mathrm{d}x
-\left[N-(p-1)a-b-1\right]
\int_{\mathbb{R}^N}
\frac{|u|^p}{|x|^{(p-1)a+b+1}}
\mathrm{d}x
\\&\quad\ge
C(N,p,b)
\inf_{c\in\mathbb{R}}
\left[\int_{\mathbb{R}^N}
\frac{\left|\nabla \left(u
-ce^{-\frac{|x|^{b-a+1}}{b-a+1}}
\right)\right|^p}
{|x|^{pb}}\mathrm{d}x
+\int_{\mathbb{R}^N}
\frac{\left|u
-ce^{-\frac{|x|^{b-a+1}}{b-a+1}}
\right|^p}
{|x|^{pa}}\mathrm{d}x\right].
\end{align*}
\end{thm}

\begin{thm}\label{thm-5}
Let $N>p\ge2$, $0\le b
<\frac{N-p}{p}$, $a<\frac{Nb}{N-p}$ and
$(p-1)a+b+1=\frac{pbN}{N-p}$.
For each $u\in\mathcal{C}_0^\infty
(\mathbb{R}^N\setminus
\{\mathbf{0}\})$, there is a universal constant $C(N,p,b)>0$ such that
\begin{align*}
&\int_{\mathbb{R}^N}
\frac{|\nabla u|^p}{|x|^{pb}}\mathrm{d}x
+(p-1)\int_{\mathbb{R}^N}
\frac{|u|^p}{|x|^{pa}}\mathrm{d}x
-\left[N-(p-1)a-b-1\right]
\int_{\mathbb{R}^N}
\frac{|u|^p}{|x|^{(p-1)a+b+1}}
\mathrm{d}x
\\&\quad
\ge
C(N,p,b)
\inf_{c\in\mathbb{R}}
\int_{\mathbb{R}^N}
\frac{\left|u
-ce^{-\frac{|x|^{b-a+1}}{b-a+1}}
\right|^p}
{|x|^{(p-1)a+b+1}}\mathrm{d}x.
\end{align*}
\end{thm}

\begin{rem}
In contrast with the previous works, some comments about Theorems \ref{thm-3}--\ref{thm-5} are listed below.
\begin{enumerate}[itemsep=0pt, topsep=0pt, parsep=0pt]

\item[(1)]
Theorem \ref{thm-3} establishes the stability of the scale non-invariant $L^p$-CKN inequalities for $1<p<2$.

\item[(2)]
If $p=2$, Theorem \ref{thm-5} reduces into \cite[Theorem 3.6]{Cazacu24}. Thereby, Theorem \ref{thm-5} improves and generalizes \cite[Theorem 3.6]{Cazacu24} from $L^2$ to $L^p$ setting.

\end{enumerate}
\end{rem}

\subsection{Main novelty, difficulties and strategy}

\noindent In this section, we will state the main novelty, difficulties and strategy of this paper.

\subsubsection{Novelty of this paper}

\begin{enumerate}
[itemsep=0pt, topsep=2pt, parsep=0pt]

\item[$(1)$]
An estimate in terms of gradient is \textbf{not} possible for the scale invariant $L^p$-CKN inequality \eqref{ckn} has been investigated in this paper. Furthermore, we also show the difference about stability results between scale invariant $L^p$-CKN inequalities and scale non-invariant $L^p$-CKN inequalities.

\item[$(2)$]
This paper establishes the stability of  $L^p$-CKN inequalities, which improve and generalize the recent papers \cite{Cazacu24,Do23}.

\item[$(3)$]
As a byproduct, we also obtain some weighted $L^p$-Poincar\'{e} inequalities (including the case $p=2$) for the log-concave measure on the Borel sets of $\mathbb{R}^N$.

\end{enumerate}

\subsubsection{Main difficulties and strategy}

\noindent Now, we briefly state main difficulties and strategy during the proof of our main results.

\emph{First of all}, the following $L^p$-CKN identities obtained in \cite[Theorem 4.2]{Do23} play a crucial role in studying the stability of $L^p$-CKN inequalities: let $N\ge1$, $p>1$, $b-a+1>0$ and $b\le\frac{N-p}{p}$, for each $u\in\mathcal{C}_0^\infty
(\mathbb{R}^N\setminus\{\mathbf{0}\})
\setminus\{0\}$,
\begin{align}\label{1.9}
&\int_{\mathbb{R}^N}
\frac{|\nabla u|^p}{|x|^{pb}}\mathrm{d}x
+(p-1)\int_{\mathbb{R}^N}
\frac{|u|^p}{|x|^{pa}}\mathrm{d}x
-\left[N-(p-1)a-b-1\right]
\int_{\mathbb{R}^N}
\frac{|u|^p}{|x|^{(p-1)a+b+1}}
\mathrm{d}x
\nonumber\\&\quad=
\int_{\mathbb{R}^N}\frac{1}
{|x|^{pb}}
\mathcal{G}_p\left(-u|x|^{b-a-1}x,
\nabla u\right)
\mathrm{d}x,
\end{align}
and
\begin{align}\label{1.10}
&\left(\int_{\mathbb{R}^N}
\frac{|\nabla u|^p}{|x|^{pb}}\mathrm{d}x
\right)^{\frac{1}{p}}
\left(\int_{\mathbb{R}^N}
\frac{|u|^p}{|x|^{pa}}
\mathrm{d}x\right)^{\frac{p-1}{p}}
-\frac{N-(p-1)a-b-1}{p}
\int_{\mathbb{R}^N}
\frac{|u|^p}{|x|^{(p-1)a+b+1}}
\mathrm{d}x
\nonumber\\&\quad=
\frac{1}{p}
\int_{\mathbb{R}^N}\frac{1}{|x|^{pb}}
\mathcal{G}_p
\left(-\left(
\frac{\int_{\mathbb{R}^N}
\frac{|\nabla u|^p}{|x|^{pb}}
\mathrm{d}x}{\int_{\mathbb{R}^N}
\frac{|u|^p}{|x|^{pa}}\mathrm{d}x}
\right)^{\frac{1}{p^2}}
u|x|^{b-a-1}x,
\left(
\frac{\int_{\mathbb{R}^N}
\frac{|u|^p}{|x|^{pa}}
\mathrm{d}x}{\int_{\mathbb{R}^N}
\frac{|\nabla u|^p}{|x|^{pb}}\mathrm{d}x}
\right)^{\frac{p-1}{p^2}}
\nabla u\right)
\mathrm{d}x,
\end{align}
where
\begin{equation*}
\mathcal{G}_p\left(\overrightarrow{X},
\overrightarrow{Y}\right)
:=\left|\overrightarrow{Y}\right|^p
-\left|\overrightarrow{X}\right|^p
-p\left|\overrightarrow{X}\right|^{p-2}\overrightarrow{X}
\cdot\left(\overrightarrow{Y}-\overrightarrow{X}\right),
\ \
\mathrm{for}
\
\mathrm{vectors}
\
\overrightarrow{X},\overrightarrow{Y}\in\mathbb{R}^N.
\end{equation*}
It deduces from \cite[Lemma 1.1]{Do23} that for all $p\ge2$, there is $C_p\in(0,1]$ such that
\[
\mathcal{G}_p\left(\overrightarrow{X},
\overrightarrow{Y}\right)
\ge C_p\left|\overrightarrow{X}
-\overrightarrow{Y}\right|^p.
\]
However, if we consider the stability of $L^p$-CKN inequalities when $1<p<2$, we need another elementary estimate of $\mathcal{G}_p$. Inspired by \cite[Lemma 2.1]{Figalli22}, for all $1<p<2$, there is $c_p>0$ such that
\begin{equation}\label{1.11}
\mathcal{G}_p\left(\overrightarrow{X},
\overrightarrow{Y}\right)
\ge c_p\min\left\{\left|\overrightarrow{Y}
-\overrightarrow{X}\right|^p,
\left|\overrightarrow{X}\right|^{p-2}
\left|\overrightarrow{Y}
-\overrightarrow{X}\right|^2\right\}.
\end{equation}

\emph{Secondly}, \eqref{1.11} is vital in proving the stability of the $L^p$-CKN inequality \eqref{ckn} for $1<p<2$. In order to use \eqref{1.11} to \eqref{1.9} and \eqref{1.10}, the first thing we need to do is to compare $\big|\overrightarrow{Y}-\overrightarrow{X}\big|$ with $\big|\overrightarrow{X}\big|$. Let us take \eqref{1.9} as an example, let $\overrightarrow{X}:=-
u|x|^{b-a-1}x$ and $\overrightarrow{Y}:=
\nabla u$. Then we divide the whole space $\mathbb{R}^N$ into two sets:
\begin{align*}
&\left\{x\in\mathbb{R}^N:\left|u|x|^{b-a-1}x\right|
\le\left|\nabla u+u|x|^{b-a-1}x\right|\right\},\\
&\left\{x\in\mathbb{R}^N:\left|u|x|^{b-a-1}x\right|
>\left|\nabla u+u|x|^{b-a-1}x\right|\right\}.
\end{align*}
Because $u\in\mathcal{C}_0^\infty
(\mathbb{R}^N\setminus\{\mathbf{0}\})\setminus\{0\}$, the above two sets are Borel sets. Thereby, we need to establish some weighted $L^p$-Poincar\'{e} inequalities for the log-concave measure on the Borel sets of $\mathbb{R}^N$. Borrowing some ideas from \cite[Lemma 4.1]{Do23}, we establish them, see Section \ref{poincare-borel} below.

\emph{Finally}, combining these new $L^p$-Poincar\'{e} inequalities with the scale invariant and the scale non-invariant $L^p$-CKN inequalities, several weighted versions of the stability for the scale invariant and the scale non-invariant $L^p$-CKN inequalities are obtained.

\subsection{Structure of this paper}

\begin{itemize}
[itemsep=0pt, topsep=0pt, parsep=0pt]

\item
In Section \ref{sect-2}, we present some crucial inequalities, including weighted Poincar\'{e} inequalities and a class of sharp vector inequalities.

\item
In Section \ref{sect-4}, we prove our main results. Theorem \ref{thm-6} will be first shown.

\begin{itemize}
[itemsep=0pt, topsep=0pt, parsep=0pt]

\item[--]
In Section \ref{sect-4.1}, we show the stability of the scale invariant $L^p$-CKN inequalities (namely, Theorems \ref{thm-1} and \ref{thm-2}).

\item[--]
In Section \ref{sect-5}, we prove the stability of the scale non-invariant $L^p$-CKN inequalities (namely, Theorems \ref{thm-3}, \ref{thm-4} and \ref{thm-5}).

\end{itemize}

\item
In Appendix \ref{sectpls}, we verify a technical inequality.

\end{itemize}

\section{Preliminary lemmas}\label{sect-2}

\noindent In this section, we establish a series of inequalities, which are vital to prove our main results. From now on, $C$ and $C(\cdot\cdot\cdot)$ are positive constants and may vary in different lines.

\subsection{Weighted Poincar\'{e} inequalities on $\mathbb{R}^N$}\label{poincare-whole}

\noindent Now, we recall some Poincar\'{e} inequalities on $\mathbb{R}^N$, which have been established in \cite{Cazacu24,Do23}.

\begin{lem}
[{\!\rm{\!\cite[Corollary 4.1]{Do23}}}]
\label{lem-2.1}
For all $\lambda>0$, $\sigma>0$, $N-p>\varrho\ge0$ and $\theta\ge\frac{N-p-\varrho}{N-p}$,
\[
\lambda^{\frac{p(N-p-\varrho)}{N-p}}
\int_{\mathbb{R}^N}\frac{|\nabla f(x)|^p}{|x|^\varrho}
e^{-\sigma\frac{|x|^\theta}
{\lambda^\theta}}\mathrm{d}x
\ge C(N,p,\theta,\sigma,\varrho)
\inf\limits_{c\in\mathbb{R}}
\int_{\mathbb{R}^N}
\frac{|f(x)-c|^p}
{|x|^{\frac{N\varrho}{N-p}}}
e^{-\sigma\frac{|x|^\theta}
{\lambda^\theta}}\mathrm{d}x.
\]
\end{lem}

Here are some examples that follows from Lemma \ref{lem-2.1} immediately by checking the appropriate pairs of $\lambda$ and $p$. By choosing $\lambda=1$ in Lemma \ref{lem-2.1}, there holds the following inequality (see also \cite[Lemma 4.1]{Do23}).

\begin{lem}\label{lem-2.2}
For all $\sigma>0$, $N-p>\varrho\ge0$ and $\theta\ge\frac{N-p-\varrho}{N-p}$,
\[
\int_{\mathbb{R}^N}\frac{|\nabla f(x)|^p}{|x|^\varrho}
e^{-\sigma|x|^\theta}\mathrm{d}x
\ge C(N,p,\theta,\sigma,\varrho)
\inf\limits_{c\in\mathbb{R}}
\int_{\mathbb{R}^N}
\frac{|f(x)-c|^p}
{|x|^{\frac{N\varrho}{N-p}}}
e^{-\sigma|x|^\theta}\mathrm{d}x.
\]
\end{lem}

Let $p=2$ in Lemma \ref{lem-2.1}, we have the following $L^2$-Poincar\'{e} inequality which also has been presented in \cite[Lemma 3.2]{Do23}.

\begin{lem}\label{lem-2.3}
For all $\lambda>0$, $\sigma>0$, $N-2>\varrho\ge0$ and $\theta\ge\frac{N-2-\varrho}{N-2}$,
\[
\lambda^{\frac{2(N-2-\varrho)}{N-2}}
\int_{\mathbb{R}^N}\frac{|\nabla f(x)|^2}{|x|^\varrho}
e^{-\sigma\frac{|x|^\theta}
{\lambda^\theta}}\mathrm{d}x
\ge C(N,\theta,\sigma,\varrho)
\inf\limits_{c\in\mathbb{R}}
\int_{\mathbb{R}^N}
\frac{|f(x)-c|^2}
{|x|^{\frac{N\varrho}{N-2}}}
e^{-\sigma
\frac{|x|^\theta}
{\lambda^\theta}}\mathrm{d}x.
\]
\end{lem}

When $p=2$ and $\lambda=1$ in Lemma \ref{lem-2.1}, the following inequality holds, also see  \cite[Lemma 3.6]{Cazacu24} or \cite[Lemma 3.1]{Do23}.

\begin{lem}\label{lem-2.4}
For all $\sigma>0$, $N-2>\varrho\ge0$ and $\theta\ge\frac{N-2-\varrho}{N-2}$,
\[
\int_{\mathbb{R}^N}\frac{|\nabla f(x)|^2}{|x|^\varrho}
e^{-\sigma|x|^\theta}\mathrm{d}x
\ge C(N,\theta,\sigma,\varrho)
\inf\limits_{c\in\mathbb{R}}
\int_{\mathbb{R}^N}
\frac{|f(x)-c|^2}
{|x|^{\frac{N\varrho}{N-2}}}
e^{-\sigma|x|^\theta}\mathrm{d}x.
\]
\end{lem}

\subsection{Weighted Poincar\'{e} inequalities on the Borel sets of $\mathbb{R}^N$}
\label{poincare-borel}

\noindent Before stating some weighted Poincar\'{e} inequalities on the Borel sets of $\mathbb{R}^N$, we first present a definition. If $U$ is a smooth convex function, then the probability measure $e^{-U}\mathrm{d}x$ on the Borel sets of $\mathbb{R}^N$ is called \emph{log-concave}, see \cite[p. 203]{Bakry14} for details.

\begin{lem}\label{lem-2.5}
Let $\Omega\subseteq\mathbb{R}^N$ be a Borel set, for all  $\sigma>0$, $N-2>\varrho\ge0$ and $\theta\ge\frac{N-2-\varrho}{N-2}$,
\[
\int_{\Omega}\frac{|\nabla f(x)|^2}{|x|^\varrho}
e^{-\sigma|x|^\theta}\mathrm{d}x
\ge C(N,\theta,\sigma,\varrho)
\inf\limits_{c\in\mathbb{R}}
\int_{\Omega}\frac{|f(x)-c|^2}
{|x|^{\frac{N\varrho}{N-2}}}
e^{-\sigma|x|^\theta}\mathrm{d}x.
\]
\end{lem}

\begin{proof}[\rm\textbf{Proof}]
For $\lambda\ge1$, let $\bar{f}(x):
=\left(\frac{1}{\lambda}
\right)^{\frac{1}{2}}
f\left(|x|^{\lambda-1}x\right)$, then
\[
|\nabla\bar{f}(x)|\le
\lambda^{\frac{1}{2}}|x|^{\lambda-1}
\left|\nabla f
\left(|x|^{\lambda-1}x\right)\right|.
\]
For convenience, denote $y:=|x|^{\lambda-1}x$, it yields from \cite[p. 6]{Lam17} that $\mathrm{d}y=\lambda|x|^{N(\lambda-1)}\mathrm{d}x$, then
\begin{align*}
\int_{\Omega}\frac{|\nabla f(y)|^2}{|y|^\varrho}
e^{-\sigma|y|^\theta}\mathrm{d}y
&=\int_{\Omega^\lambda}
\frac{\left|\nabla f
\left(|x|^{\lambda-1}x\right)
\right|^2}
{|x|^{\lambda\varrho}}
e^{-\sigma|x|^{\lambda\theta}}
\lambda|x|^{N(\lambda-1)}\mathrm{d}x
\\&\ge\int_{\Omega^\lambda}
\frac{
\left|\nabla\bar{f}(x)\right|^2}
{\lambda|x|^{2(\lambda-1)
+\lambda\varrho-N(\lambda-1)}}
e^{-\sigma|x|^{\lambda\theta}}
\lambda\mathrm{d}x
\\&=\int_{\Omega^\lambda}
\frac{
\left|\nabla\bar{f}(x)\right|^2}
{|x|^{\lambda(2+\varrho-N)+N-2}}
e^{-\sigma|x|^{\lambda\theta}}
\mathrm{d}x,
\end{align*}
where $\Omega^\lambda:
=\{x:|x|^{\lambda-1}x\in\Omega\}$ is a Borel set. Choosing $\lambda=\frac{N-2}{N-2-\varrho}\ge 1$, it is easy to check that the measure $\exp\left[-\sigma|x|^{\frac{N-2}
{N-2-\varrho}\theta}\right]\mathrm{d}x$ is log-concave for all $\theta\ge\frac{N-2-\varrho}{N-2}$ and $\sigma>0$. Then, with the aid of \cite[Theorem 4.6.3]{Bakry14}, it gives that
\begin{align*}
\int_{\Omega}\frac{|\nabla f(y)|^2}{|y|^\varrho}
e^{-\sigma|y|^\theta}\mathrm{d}y
&\ge\int_{\Omega^\lambda}
|\nabla\bar{f}(x)|^2
e^{-\sigma|x|^{\frac{N-2}
{N-2-\varrho}\theta}}\mathrm{d}x
\\&\ge C(N,\sigma,\theta,\varrho)
\inf_{c\in\mathbb{R}}
\int_{\Omega^\lambda}
|\bar{f}(x)-c|^2
e^{-\sigma|x|^{\lambda\theta}}
\mathrm{d}x
\\&= C(N,\sigma,\theta,\varrho)
\inf_{c\in\mathbb{R}}
\int_{\Omega^\lambda}
\frac{\left|f
\left(|x|^{\lambda-1}x\right)
-c\right|^2}
{|x|^{N(\lambda-1)}}
e^{-\sigma|x|^{\lambda\theta}}
|x|^{N(\lambda-1)}
\mathrm{d}x
\\&= C(N,\sigma,\theta,\varrho)
\inf_{c\in\mathbb{R}}
\int_{\Omega}
\frac{\left|f\left(y\right)
-c\right|^2}
{|y|^{\frac{N\varrho}{N-2}}}
e^{-\sigma|y|^{\theta}}\mathrm{d}y,
\end{align*}
which are our desired estimates. Hence, the proof is completed.
\end{proof}

We obtain the following inequalities by applying the scaling argument to Lemma \ref{lem-2.5}.

\begin{lem}\label{lem-2.6}
Let $\Omega\subseteq\mathbb{R}^N$ be a Borel set, for all  $\lambda>0$, $\sigma>0$, $N-2>\varrho\ge0$ and $\theta\ge\frac{N-2-\varrho}{N-2}$,
\[
\lambda^{\frac{2(N-2-\varrho)}{N-2}}
\int_{\Omega}\frac{|\nabla f(x)|^2}{|x|^\varrho}
e^{-\sigma\frac{|x|^\theta}
{\lambda^\theta}}\mathrm{d}x
\ge C(N,\theta,\sigma,\varrho)
\inf\limits_{c\in\mathbb{R}}
\int_{\Omega}\frac{|f(x)-c|^2}
{|x|^{\frac{N\varrho}{N-2}}}
e^{-\sigma\frac{|x|^\theta}
{\lambda^\theta}}\mathrm{d}x.
\]
\end{lem}

The following $L^p$-Poincar\'{e} inequalities on the  Borel sets of $\mathbb{R}^N$ are also required to research the stability of $L^p$-CKN inequalities. Similar arguments as those of Lemma \ref{lem-2.5} or \cite[Lemma 4.1]{Do23}, we obtain them.

\begin{lem}\label{lem-2.7}
Let $\Omega\subseteq\mathbb{R}^N$ be a Borel set, for all $\sigma>0$, $N-p>\varrho\ge0$ and $\theta\ge\frac{N-p-\varrho}{N-p}$,
\[
\int_{\Omega}\frac{|\nabla f(x)|^p}{|x|^\varrho}
e^{-\sigma|x|^\theta}\mathrm{d}x
\ge C(N,p,\theta,\sigma,\varrho)
\inf\limits_{c\in\mathbb{R}}
\int_{\Omega}\frac{|f(x)-c|^p}
{|x|^{\frac{N\varrho}{N-p}}}
e^{-\sigma|x|^\theta}\mathrm{d}x.
\]
\end{lem}

Similarly, the following weighted $L^p$-Poincar\'{e} inequalities hold by the scaling argument.

\begin{lem}\label{lem-2.8}
Let $\Omega\subseteq\mathbb{R}^N$ be a Borel set, for all $\lambda>0$, $\sigma>0$, $N-p>\varrho\ge0$ and $\theta\ge\frac{N-p-\varrho}{N-p}$,
\[
\lambda^{\frac{p(N-p-\varrho)}{N-p}}
\int_{\Omega}\frac{|\nabla f(x)|^p}{|x|^\varrho}
e^{-\sigma\frac{|x|^\theta}
{\lambda^\theta}}\mathrm{d}x
\ge C(N,p,\theta,\sigma,\varrho)
\inf\limits_{c\in\mathbb{R}}
\int_{\Omega}\frac{|f(x)-c|^p}
{|x|^{\frac{N\varrho}{N-p}}}
e^{-\sigma\frac{|x|^\theta}
{\lambda^\theta}}\mathrm{d}x.
\]
\end{lem}

\subsection{Sharp vector inequalities}
\label{sect-2.3}

\noindent In order to study the stability of $L^p$-CKN inequalities for $1<p<2$, we borrow some ideas from the following lemma that proposed by Figalli and Zhang in \cite[Lemma 2.1]{Figalli22}.

\begin{lem}
[{\!\rm{\!\cite[Lemma 2.1]{Figalli22}}}]
\label{lem-2.9}
Assume that $x,y\in\mathbb{R}^N$. Then, for each $\gamma>0$, there is $c_{p,\gamma}>0$ such that
\begin{enumerate}[itemsep=0pt, topsep=1pt, parsep=0pt]

\item[(1)] if $1<p<2$,
\begin{align*}
|x+y|^p
-|x|^p
-p|x|^{p-2}x\cdot y
&\ge\frac{1-\gamma}{2}
\left[p|x|^{p-2}|y|^2
+p(p-2)|z|^{p-2}
\left(|x|-|x+y|\right)^2\right]
\\&\quad+c_{p,\gamma}
\min\left\{|y|^p,|x|^{p-2}
|y|^2\right\},
\end{align*}
where
\[
z=z(x,x+y)
:=\begin{cases}
\left[\frac{|x+y|}
{(2-p)|x+y|+(p-1)|x|}
\right]^{\frac{1}{p-2}}x,
\ \ &if\ |x|<|x+y|;\\
x,
\ \ &if\ |x|\ge|x+y|;
\end{cases}
\]

\item[(2)] if $p\ge 2$,
\begin{align*}
|x+y|^p
-|x|^p
-p|x|^{p-2}x\cdot y
&\ge\frac{1-\gamma}{2}
\left[p|x|^{p-2}|y|^2
+p(p-2)|z|^{p-2}
\left(|x|-|x+y|\right)^2\right]
\\&\quad
+c_{p,\gamma}|y|^p,
\end{align*}
where
\[
z=z(x,x+y)
:=\begin{cases}
x,
\ \ &if\ |x|<|x+y|;\\
\left(\frac{|x+y|}{|x|}
\right)^{\frac{1}{p-2}}(x+y),
\ \ &if\ |x|\ge|x+y|.
\end{cases}
\]
\end{enumerate}
\end{lem}

Due to the arbitrariness of $\gamma$, let us choose $\gamma=1$ in Lemma \ref{lem-2.9}, we get the following result.

\begin{cor}\label{cor-2.10}
Assume that $x,y\in\mathbb{R}^N$. Then there is $c_p>0$ such that
\begin{equation*}
|x+y|^p
-|x|^p
-p|x|^{p-2}x\cdot y\ge
\begin{cases}
c_p\min\left\{|y|^p,
|x|^{p-2}|y|^2\right\},
\ \ &if\ 1<p<2;\\[1mm]
c_p|y|^p,
\ \ &if\ p\ge 2.
\end{cases}
\end{equation*}
\end{cor}

Given $N\ge 1$, $p>1$, $\overrightarrow{X}$ and $\overrightarrow{Y}$ be vectors on $\mathbb{R}^N$, and let
\begin{align}\label{defRp}
\mathcal{G}_p\left(\overrightarrow{X},
\overrightarrow{Y}\right)
:=\left|\overrightarrow{Y}\right|^p
-\left|\overrightarrow{X}\right|^p
-p\left|\overrightarrow{X}\right|^{p-2}
\overrightarrow{X}
\cdot\left(\overrightarrow{Y}
-\overrightarrow{X}\right).
\end{align}
Let $x+y=\overrightarrow{Y}$ and $x=\overrightarrow{X}$ in Corollary \ref{cor-2.10}, we see that
\begin{equation}\label{2.2}
\mathcal{G}_p\left(\overrightarrow{X},
\overrightarrow{Y}\right)\ge
\begin{cases}
c_p\min\left\{\left|\overrightarrow{Y}
-\overrightarrow{X}\right|^p,
\left|\overrightarrow{X}\right|^{p-2}
\left|\overrightarrow{Y}
-\overrightarrow{X}\right|^2\right\},
\ \ &\mathrm{if}\ 1<p<2;\\[1mm]
c_p\left|\overrightarrow{Y}-\overrightarrow{X}\right|^p,
\ \ &\mathrm{if}\ p\ge 2,
\end{cases}
\end{equation}
for some $c_p>0$.

\section{Proof of main results}\label{sect-4}

\noindent In this section, we prove our main results. We first prove Theorem \ref{thm-6}.

\begin{proof}[\rm\textbf{Proof of Theorem \ref{thm-6}}]
For any $u\in\mathcal{C}_0^\infty
(\mathbb{R}^N\setminus
\{\mathbf{0}\})$ and $\lambda>0$, we define
\begin{align*}
u_\lambda(x):=\lambda^{\frac{N-(p-1)a-b-1}{p}}
u(\lambda x).
\end{align*}
By direct calculations, we have
\begin{align*}
\int_{\mathbb{R}^N}
\frac{|u_\lambda|^p}{|x|^{(p-1)a+b+1}}
\mathrm{d}x
&=\int_{\mathbb{R}^N}
\frac{|u|^p}{|x|^{(p-1)a+b+1}}
\mathrm{d}x,
\\
\int_{\mathbb{R}^N}
\frac{|\nabla u_\lambda|^p}{|x|^{pb}}\mathrm{d}x
&=\lambda^{(p-1)(b+1-a)}
\int_{\mathbb{R}^N}
\frac{|\nabla u|^p}{|x|^{pb}}\mathrm{d}x,
\\
\int_{\mathbb{R}^N}
\frac{|u_\lambda|^p}{|x|^{pa}}
\mathrm{d}x
&=\lambda^{a-b-1}\int_{\mathbb{R}^N}
\frac{|u|^p}{|x|^{pa}}
\mathrm{d}x.
\end{align*}
Therefore, the $L^p$-CKN inequality \eqref{ckn} is scale invariant, that is,
\begin{align*}
\delta(u_\lambda)=\delta(u),
\quad \mbox{for any}\ u\in\mathcal{C}_0^\infty
(\mathbb{R}^N\setminus
\{\mathbf{0}\}),\ \lambda>0,
\end{align*}
where
\begin{align*}
\delta(u):=\left(\int_{\mathbb{R}^N}
\frac{|\nabla u|^p}{|x|^{pb}}\mathrm{d}x
\right)^{\frac{1}{p}}
\left(\int_{\mathbb{R}^N}
\frac{|u|^p}{|x|^{pa}}
\mathrm{d}x\right)^{\frac{p-1}{p}}
-\frac{N-(p-1)a-b-1}{p}
\int_{\mathbb{R}^N}
\frac{|u|^p}{|x|^{(p-1)a+b+1}}
\mathrm{d}x.
\end{align*}

Now, let $C_1, C_2$ be given. Suppose $C_1>0$. Let $u\in\mathcal{C}_0^\infty
(\mathbb{R}^N\setminus
\{\mathbf{0}\})\setminus\mathcal{M}_{a,b}$ satisfy
\begin{align*}
\int_{\mathbb{R}^N}
\frac{|\nabla u|^p}{|x|^{pb}}\mathrm{d}x
\leq 2\inf_{v\in\mathcal{M}_{a,b}}
\int_{\mathbb{R}^N}
\frac{|\nabla (u-v)|^p}{|x|^{pb}}\mathrm{d}x.
\end{align*}
By the scaling invariance of $\delta$, there always exists a sufficiently small $\lambda$ depending upon $C_1$ and $\inf_{v\in\mathcal{M}_{a,b}}\int_{\mathbb{R}^N}
\frac{|\nabla (u-v)|^p}{|x|^{pb}}\mathrm{d}x$ such that the following holds,
\begin{align*}
\delta(u)=\delta(u_\lambda)
&\leq \frac{C_1}{2}\lambda^{-(p-1)(b+1-a)}\int_{\mathbb{R}^N}
\frac{|\nabla u_\lambda|^p}{|x|^{pb}}\mathrm{d}x
\\
& =\frac{C_1}{2}
\int_{\mathbb{R}^N}
\frac{|\nabla u|^p}{|x|^{pb}}\mathrm{d}x
\leq C_1\inf_{v\in\mathcal{M}_{a,b}}
\int_{\mathbb{R}^N}
\frac{|\nabla (u-v)|^p}{|x|^{pb}}\mathrm{d}x.
\end{align*}

Next, suppose $C_1>0$, and let $u\in\mathcal{C}_0^\infty
(\mathbb{R}^N\setminus
\{\mathbf{0}\})\setminus\mathcal{M}_{a,b}$ satisfy
\begin{align*}
\int_{\mathbb{R}^N}
\frac{
\left|u\right|^p}{|x|^{pa}}
\mathrm{d}x
\leq 2\inf_{v\in\mathcal{M}_{a,b}}
\int_{\mathbb{R}^N}
\frac{
\left|u-v\right|^p}{|x|^{pa}}
\mathrm{d}x.
\end{align*}
A symmetric calculation for sufficiently large $\lambda$ such that the following also holds,
\begin{align*}
\delta(u)=\delta(u_\lambda)
\leq \frac{C_2}{2}\lambda^{b+1-a}\int_{\mathbb{R}^N}
\frac{
\left|u_\lambda\right|^p}{|x|^{pa}}
\mathrm{d}x
=\frac{C_2}{2}\int_{\mathbb{R}^N}
\frac{
\left|u\right|^p}{|x|^{pa}}
\mathrm{d}x
\leq C_2
\inf_{v\in\mathcal{M}_{a,b}}
\int_{\mathbb{R}^N}
\frac{\left|u-v\right|^p}{|x|^{pa}}
\mathrm{d}x.
\end{align*}
The proof is thereby completed.
\end{proof}

In the following of this section, we will use the $L^p$-CKN identities derived in Lemma \ref{lem-4.1} below and the Poincar\'{e} inequality in Section \ref{sect-2} to study the stability of $L^p$-CKN inequalities. This idea was developed by Cazacu \emph{et al.} to derive the stability of CKN inequalities (see \cite{Cazacu24,Do23}). Before we present the proof, we need the following $L^p$-CKN identities.

\begin{lem}
[{\!\rm{\!\cite[Theorem 4.2]{Do23}}}]
\label{lem-4.1}
Let $N\ge1$, $p>1$, $b-a+1>0$ and $b\le\frac{N-p}{p}$. For each $u\in\mathcal{C}_0^\infty
(\mathbb{R}^N\setminus\{\mathbf{0}\})
\setminus\{0\}$, then
\begin{align}\label{4.1}
&\int_{\mathbb{R}^N}
\frac{|\nabla u|^p}{|x|^{pb}}\mathrm{d}x
+(p-1)\int_{\mathbb{R}^N}
\frac{|u|^p}{|x|^{pa}}\mathrm{d}x
-\left[N-(p-1)a-b-1\right]
\int_{\mathbb{R}^N}
\frac{|u|^p}{|x|^{(p-1)a+b+1}}
\mathrm{d}x
\nonumber\\&\quad=
\int_{\mathbb{R}^N}\frac{1}
{|x|^{pb}}
\mathcal{G}_p\left(-u|x|^{b-a-1}x,
\nabla u\right)
\mathrm{d}x,
\end{align}
and
\begin{align}\label{4.2}
&\left(\int_{\mathbb{R}^N}
\frac{|\nabla u|^p}{|x|^{pb}}\mathrm{d}x
\right)^{\frac{1}{p}}
\left(\int_{\mathbb{R}^N}
\frac{|u|^p}{|x|^{pa}}
\mathrm{d}x\right)^{\frac{p-1}{p}}
-\frac{N-(p-1)a-b-1}{p}
\int_{\mathbb{R}^N}
\frac{|u|^p}{|x|^{(p-1)a+b+1}}
\mathrm{d}x
\nonumber\\&\quad=
\frac{1}{p}
\int_{\mathbb{R}^N}\frac{1}{|x|^{pb}}
\mathcal{G}_p\left(-\left(
\frac{\int_{\mathbb{R}^N}
\frac{|\nabla u|^p}{|x|^{pb}}\mathrm{d}x}
{\int_{\mathbb{R}^N}
\frac{|u|^p}{|x|^{pa}}
\mathrm{d}x}
\right)^{\frac{1}{p^2}}u|x|^{b-a-1}x,
\left(
\frac{\int_{\mathbb{R}^N}
\frac{|u|^p}{|x|^{pa}}
\mathrm{d}x}{\int_{\mathbb{R}^N}
\frac{|\nabla u|^p}{|x|^{pb}}\mathrm{d}x}
\right)^{\frac{p-1}{p^2}}\nabla u\right)
\mathrm{d}x,
\end{align}
where $\mathcal{G}_p$ is defined in \eqref{defRp}.
\end{lem}

\begin{rem}
Because the assumptions of Theorems \ref{thm-1}, \ref{thm-2}, \ref{thm-3}, \ref{thm-4} and \ref{thm-5} all satisfy those of Lemma \ref{lem-4.1}, then Lemma \ref{lem-4.1} can be utilized to prove these theorems.
\end{rem}

\subsection{The stability results of the scale invariant $L^p$-CKN inequalities: proof of Theorems \ref{thm-1} and \ref{thm-2}}\label{sect-4.1}

\noindent The subsection is concerned about the stability results of the scale invariant $L^p$-CKN inequalities. We begin with proving Theorem \ref{thm-1}.

\begin{proof}[\rm\textbf{Proof of Theorem \ref{thm-1}}]
For simplicity, let
\begin{equation*}
\lambda=\left(
\frac{\int_{\mathbb{R}^N}
\frac{|u|^p}{|x|^{pa}}
\mathrm{d}x}{\int_{\mathbb{R}^N}
\frac{|\nabla u|^p}{|x|^{pb}}\mathrm{d}x}
\right)^{\frac{1}{p(b-a+1)}}>0.
\end{equation*}
From \eqref{4.2}, we obtain
\begin{align}\label{4.3}
&\bigg(\int_{\mathbb{R}^N}
\frac{|\nabla u|^p}{|x|^{pb}}\mathrm{d}x
\bigg)^{\frac{1}{p}}
\left(\int_{\mathbb{R}^N}
\frac{|u|^p}{|x|^{pa}}
\mathrm{d}x\right)^{\frac{p-1}{p}}
-\frac{N-(p-1)a-b-1}{p}
\int_{\mathbb{R}^N}
\frac{|u|^p}{|x|^{(p-1)a+b+1}}
\mathrm{d}x
\nonumber\\&\quad=
\frac{1}{p}
\int_{\mathbb{R}^N}\frac{1}{|x|^{pb}}
\mathcal{G}_p
\left(-\lambda^{-\frac{b-a+1}{p}}
u|x|^{b-a-1}x,
\lambda^{\frac{(p-1)(b-a+1)}{p}}
\nabla u\right)
\mathrm{d}x
\nonumber\\&\quad=
\frac{1}{p}
\int_{\Omega_{u,1}}\frac{1}
{|x|^{pb}}
\mathcal{G}_p
\left(-\lambda^{-\frac{b-a+1}{p}}
u|x|^{b-a-1}x,
\lambda^{\frac{(p-1)(b-a+1)}{p}}
\nabla u\right)
\mathrm{d}x
\nonumber\\&\qquad
+\frac{1}{p}
\int_{\Omega_{u,2}}\frac{1}{|x|^{pb}}
\mathcal{G}_p
\left(-\lambda^{-\frac{b-a+1}{p}}
u|x|^{b-a-1}x,
\lambda^{\frac{(p-1)(b-a+1)}{p}}
\nabla u\right)
\mathrm{d}x
\nonumber\\&\quad
=\frac{1}{p}\Pi_1+\frac{1}{p}\Pi_2,
\end{align}
where
\begin{align*}
\Omega_{u,1}
&:=\left\{x\in\mathbb{R}^N:
\left|u|x|^{b-a-1}x\right|
\le\left|\lambda^{b-a+1}\nabla u+u|x|^{b-a-1}x\right|
\right\},\\
\Omega_{u,2}&
:=\left\{x\in\mathbb{R}^N:
\left|u|x|^{b-a-1}x\right|
>\left|\lambda^{b-a+1}\nabla u+u|x|^{b-a-1}x\right|
\right\},
\end{align*}
and $\Omega_{u,1}\cup\Omega_{u,2}
=\mathbb{R}^N$. Since $u\in\mathcal{C}_0^\infty
(\mathbb{R}^N
\setminus\{\mathbf{0}\})
\setminus\{0\}$, $\Omega_{u,1}$ and $\Omega_{u,2}$ are Borel sets, and then we can use Lemmas \ref{lem-2.6} and \ref{lem-2.8} to estimate $\Pi_1$ and $\Pi_2$.

For $\Pi_1$. We deduce from \eqref{2.2} and $\left|u|x|^{b-a-1}x\right|
\le\left|\lambda^{b-a+1}\nabla u+u|x|^{b-a-1}x\right|$ that
\begin{align*}
&\mathcal{G}_p
\left(-\lambda^{-\frac{b-a+1}{p}}
u|x|^{b-a-1}x,
\lambda^{\frac{(p-1)(b-a+1)}{p}}
\nabla u\right)
\nonumber\\&\quad
\ge c_p\lambda^{(p-1)(b-a+1)}
\left|\nabla u
+\lambda^{-(b-a+1)}
u|x|^{b-a-1}x\right|^p.
\end{align*}
This together with \eqref{4.3} and Lemma \ref{lem-2.8} (with $\varrho=pb$, $\sigma=\frac{p}{b-a+1}$ and $\theta=b-a+1$), we obtain
\begin{align}\label{4.4}
\Pi_1&\ge c_p\lambda^{(p-1)(b-a+1)}
\int_{\Omega_{u,1}}
\frac{1}{|x|^{pb}}
\left|\nabla u
+\lambda^{-(b-a+1)}
u|x|^{b-a-1}x\right|^p
\mathrm{d}x
\nonumber\\&
=c_p\lambda^{(p-1)(b-a+1)}
\int_{\Omega_{u,1}}
\frac{1}{|x|^{pb}}
\left|\nabla
\left(ue^{\frac{|x|^{b-a+1}}
{(b-a+1)\lambda^{b-a+1}}}\right)
\right|^p
e^{{-\frac{p|x|^{b-a+1}}
{(b-a+1)\lambda^{b-a+1}}}}
\mathrm{d}x
\nonumber\\&\ge
C(N,p,b)\lambda^{(p-1)(b-a+1)
-\left(p-\frac{p^2b}{N-p}\right)}
\inf_{c\in\mathbb{R}}
\int_{\Omega_{u,1}}
\frac{
\left|ue^{\frac{|x|^{b-a+1}}
{(b-a+1)\lambda^{b-a+1}}}-c
\right|^p}{|x|^{\frac{pbN}{N-p}}}
e^{{-\frac{p|x|^{b-a+1}}{(b-a+1)
\lambda^{b-a+1}}}}
\mathrm{d}x
\nonumber\\&=
C(N,p,b)\lambda^{(p-1)(b-a+1)
-\left(p-\frac{p^2b}{N-p}\right)}
\inf_{c\in\mathbb{R}}
\int_{\Omega_{u,1}}
\frac{
\left|u-ce^{-\frac{|x|^{b-a+1}}
{(b-a+1)\lambda^{b-a+1}}}
\right|^p}{|x|^{\frac{pbN}{N-p}}}
\mathrm{d}x.
\end{align}

For $\Pi_2$. By \eqref{2.2} and $\left|u|x|^{b-a-1}x\right|
>\left|\lambda^{b-a+1}\nabla u+u|x|^{b-a-1}x\right|$, we deduce that
\begin{align*}
&\mathcal{G}_p
\left(-\lambda^{-\frac{b-a+1}{p}}
u|x|^{b-a-1}x,
\lambda^{\frac{(p-1)(b-a+1)}{p}}
\nabla u\right)
\nonumber\\&\quad
\ge c_p\lambda^{-\frac{(p-2)(b-a+1)}{p}}
|u|^{p-2}|x|^{(p-2)(b-a)}
\left[
\lambda^{\frac{(p-1)(b-a+1)}{p}}
\nabla u
+\lambda^{-\frac{b-a+1}{p}}
u|x|^{b-a-1}x\right]^2
\nonumber\\&\quad=
c_p\lambda^{b-a+1}
|u|^{p-2}|x|^{(p-2)(b-a)}
\left[\nabla u
+\lambda^{-(b-a+1)}
u|x|^{b-a-1}x\right]^2.
\end{align*}
Then, from \eqref{4.3} and Lemma \ref{lem-2.6} (with $\varrho=pa+2b-2a$, $\theta=b-a+1$ and $\sigma=\frac{p}{b-a+1}$), we get
\begin{small}
\begin{align}\label{4.5}
\Pi_2&
\ge c_p\lambda^{b-a+1}
\int_{\Omega_{u,2}}
\frac{1}{|x|^{pb}}
|u|^{p-2}|x|^{(p-2)(b-a)}
\left[\nabla u
+\lambda^{-(b-a+1)}
u|x|^{b-a-1}x\right]^2
\mathrm{d}x
\nonumber\\&=
c_p\lambda^{b-a+1}
\int_{\Omega_{u,2}}
\frac{1}{|x|^{pa+2b-2a}}
\left[\nabla u\cdot |u|^{\frac{p-2}{2}}
+\lambda^{-(b-a+1)}
|u|^{\frac{p-2}{2}}u
|x|^{b-a-1}x\right]^2
\mathrm{d}x
\nonumber\\&=
c_p\left(\frac{2}{p}\right)^2
\lambda^{b-a+1}
\int_{\Omega_{u,2}}
\frac{1}{|x|^{pa+2b-2a}}
\left|\nabla
\left(|u|^{\frac{p-2}{2}}u
e^{\frac{\frac{p}{2}|x|^{b-a+1}}
{(b-a+1)\lambda^{b-a+1}}}\right)
\right|^2
e^{{-\frac{p|x|^{b-a+1}}
{(b-a+1)\lambda^{b-a+1}}}}
\mathrm{d}x
\nonumber\\&\ge
C(N,p,b)\lambda^{(b-a+1)
-\left(2
-\frac{2pa+4b-4a}{N-2}\right)}
\inf_{c\in\mathbb{R}}
\int_{\Omega_{u,2}}
\frac{
\left||u|^{\frac{p-2}{2}}u
e^{\frac{\frac{p}{2}|x|^{b-a+1}}
{(b-a+1)\lambda^{b-a+1}}}
-|c|^{\frac{p-2}{2}}c
\right|^2}{|x|^{\frac{N(pa+2b-2a)}
{N-2}}}
e^{{-\frac{p|x|^{b-a+1}}
{(b-a+1)\lambda^{b-a+1}}}}
\mathrm{d}x
\nonumber\\&=
C(N,p,b)\lambda^{(b-a+1)
-\left(2
-\frac{2pa+4b-4a}{N-2}\right)}
\inf_{c\in\mathbb{R}}
\int_{\Omega_{u,2}}
\frac{
\left||u|^{\frac{p-2}{2}}u
-|c|^{\frac{p-2}{2}}c
e^{-\frac{\frac{p}{2}
|x|^{b-a+1}}{(b-a+1)\lambda^{b-a+1}}}
\right|^2}{|x|^{\frac{N(pa+2b-2a)}{N-2}}}
\mathrm{d}x.
\end{align}
\end{small}

For $1<p<2$ and $m,n\in\mathbb{R}$,
\begin{equation}\label{4.6}
\left||m|^{\frac{p-2}{2}}m
-|n|^{\frac{p-2}{2}}n\right|^2
\le 4|m-n|^p.
\end{equation}
The detailed proof of \eqref{4.6} is shown in Lemma \ref{lem-a}. Therefore, substituting \eqref{4.4} and \eqref{4.5} into \eqref{4.3}, we get
\begin{align}\label{4.7}
&\left(\int_{\mathbb{R}^N}
\frac{|\nabla u|^p}{|x|^{pb}}\mathrm{d}x
\right)^{\frac{1}{p}}
\left(\int_{\mathbb{R}^N}
\frac{|u|^p}{|x|^{pa}}
\mathrm{d}x\right)^{\frac{p-1}{p}}
-\frac{N-(p-1)a-b-1}{p}
\int_{\mathbb{R}^N}
\frac{|u|^p}{|x|^{(p-1)a+b+1}}
\mathrm{d}x
\nonumber\\&\quad\ge
C(N,p,b)\lambda^{(p-1)(b-a+1)
-\left(p-\frac{p^2b}{N-p}\right)}
\inf_{c\in\mathbb{R}}
\int_{\Omega_{u,1}}
\frac{
\left|u-ce^{-\frac{|x|^{b-a+1}}
{(b-a+1)\lambda^{b-a+1}}}
\right|^p}{|x|^{\frac{pbN}{N-p}}}
\mathrm{d}x
\nonumber\\&\qquad
+C(N,p,b)\lambda^{(b-a+1)
-\left(2
-\frac{2pa+4b-4a}{N-2}\right)}
\inf_{c\in\mathbb{R}}
\int_{\Omega_{u,2}}
\frac{
\left||u|^{\frac{p-2}{2}}u
-|c|^{\frac{p-2}{2}}c
e^{-\frac{\frac{p}{2}
|x|^{b-a+1}}{(b-a+1)\lambda^{b-a+1}}}
\right|^2}{|x|^{\frac{N(pa+2b-2a)}{N-2}}}
\mathrm{d}x
\nonumber\\&\quad=
C(N,p,b)\lambda^{-(b-a+1)}
\inf_{c\in\mathbb{R}}
\int_{\Omega_{u,1}}
\frac{
\left|u-ce^{-\frac{|x|^{b-a+1}}
{(b-a+1)\lambda^{b-a+1}}}
\right|^p}{|x|^{pa}}
\mathrm{d}x
\nonumber\\&\qquad
+C(N,p,b)\lambda^{-(b-a+1)}
\inf_{c\in\mathbb{R}}
\int_{\Omega_{u,2}}
\frac{
\left||u|^{\frac{p-2}{2}}u
-|c|^{\frac{p-2}{2}}c
e^{-\frac{\frac{p}{2}
|x|^{b-a+1}}{(b-a+1)\lambda^{b-a+1}}}
\right|^2}{|x|^{pa}}
\mathrm{d}x
\nonumber\\&\quad\ge
C(N,p,b)\left(
\frac{\int_{\mathbb{R}^N}
\frac{|\nabla u|^p}{|x|^{pb}}\mathrm{d}x}
{\int_{\mathbb{R}^N}
\frac{|u|^p}{|x|^{pa}}
\mathrm{d}x}
\right)^{\frac{1}{p}}
\inf_{c\in\mathbb{R},\lambda>0}
\int_{\mathbb{R}^N}
\frac{
\left||u|^{\frac{p-2}{2}}u
-|c|^{\frac{p-2}{2}}c
e^{-\frac{\frac{p}{2}
|x|^{b-a+1}}{(b-a+1)
\lambda^{b-a+1}}}
\right|^2}{|x|^{pa}}
\mathrm{d}x.
\end{align}
The desired result follows.
\end{proof}

\begin{rem}\label{rem-1}
Under the assumption $1<p<2$, we analyze the source of other assumptions of Theorem \ref{thm-1}.

\begin{enumerate}[itemsep=0pt, topsep=1pt, parsep=0pt]

\item[$(1)$] In order to apply Lemma \ref{lem-4.1} to  obtain \eqref{4.3}, we need
\begin{equation}\label{4.8}
N\ge1,
\
p>1,
\
b-a+1>0,
\
b\le\frac{N-p}{p}.
\end{equation}

\item[$(2)$] For $\Pi_1$, in view of Lemma \ref{lem-2.8} (with $\varrho=pb$, $\sigma=\frac{p}{b-a+1}$ and $\theta=b-a+1$), there hold
\begin{small}
\begin{equation}\label{4.9}
\begin{cases}
\sigma=\frac{p}{b-a+1}>0,\\
N-p>\varrho=pb\ge 0,\\
\theta=b-a+1\ge\frac{N-p-\varrho}{N-p}
=1-\frac{pb}{N-p},
\end{cases}
\Longleftrightarrow
N>p,
\
0\le b<\frac{N-p}{p},
\
a\le\frac{Nb}{N-p}.
\end{equation}
\end{small}

\item[$(3)$] For $\Pi_2$, based on Lemma \ref{lem-2.6} (with $\varrho=pa+2b-2a$, $\theta=b-a+1$ and $\sigma=\frac{p}{b-a+1}$),
\begin{equation}\label{4.10}
\begin{cases}
\sigma=\frac{p}{b-a+1}>0,\\
N-2>\varrho=pa+2b-2a\ge 0,\\
\theta=b-a+1\ge\frac{N-2-\varrho}{N-2}
=1-\frac{pa+2b-2a}{N-2},
\end{cases}
\Longleftrightarrow
\begin{cases}
N>2,\\
N-2>pa+2b-2a\ge 0,\\
a\le\frac{Nb}{N-p},
\end{cases}
\end{equation}

\item[$(4)$] In order to obtain \eqref{4.7}, we need
\begin{equation}\label{4.11}
\frac{pb}{N-p}=\frac{pa+2b-2a}{N-2}
\
\stackrel{1<p<2}{\Longleftrightarrow}
\
a=\frac{Nb}{N-p}.
\end{equation}

\item[$(5)$] Taking the intersection of the sets proposed in \eqref{4.8}$-$\eqref{4.11} and $1<p<2$,
\begin{equation}\label{4.12}
N>2,\
1<p<2,\
0\le b
<\frac{N-p}{p},\
a=\frac{Nb}{N-p}.
\end{equation}

\item[($6$)] There are many parameters satisfying all the  assumptions of Theorem \ref{thm-1} (that is, \eqref{4.12}), such as $b=\frac{N-2}{p}$ and $a=\frac{N(N-2)}{p(N-p)}$ for $N>2,\ 1<p<2$.
\end{enumerate}
\end{rem}

We now proceed to prove Theorem \ref{thm-2} with the help of Lemmas \ref{lem-2.1} and \ref{lem-4.1}.

\begin{proof}
[\rm\textbf{Proof of Theorem \ref{thm-2}}]
Lemma \ref{lem-2.1} (with $\varrho=pb$,  $\sigma=\frac{p}{b-a+1}$ and $\theta=b-a+1$), \eqref{2.2} and \eqref{4.2} imply that
\begin{align*}
&\bigg(\int_{\mathbb{R}^N}
\frac{|\nabla u|^p}{|x|^{pb}}\mathrm{d}x
\bigg)^{\frac{1}{p}}
\left(\int_{\mathbb{R}^N}
\frac{|u|^p}{|x|^{pa}}
\mathrm{d}x\right)^{\frac{p-1}{p}}
-\frac{N-(p-1)a-b-1}{p}
\int_{\mathbb{R}^N}
\frac{|u|^p}{|x|^{(p-1)a+b+1}}
\mathrm{d}x
\nonumber\\&\quad
\ge \frac{c_p}{p}\lambda^{(p-1)(b-a+1)}
\int_{\mathbb{R}^N}\frac{1}{|x|^{pb}}
\left|\nabla u
+\lambda^{-(b-a+1)}
u|x|^{b-a-1}x\right|^p
\mathrm{d}x
\nonumber\\&\quad
=\frac{c_p}{p}\lambda^{(p-1)(b-a+1)}
\int_{\mathbb{R}^N}
\frac{1}{|x|^{pb}}
\left|\nabla
\left(ue^{\frac{|x|^{b-a+1}}
{(b-a+1)\lambda^{b-a+1}}}\right)
\right|^p
e^{{-\frac{p|x|^{b-a+1}}
{(b-a+1)\lambda^{b-a+1}}}}
\mathrm{d}x
\nonumber\\&\quad
\ge
C(N,p,b)\lambda^{(p-1)(b-a+1)
-\left(p-\frac{p^2b}{N-p}\right)}
\inf_{c\in\mathbb{R}}
\int_{\mathbb{R}^N}
\frac{
\left|ue^{\frac{|x|^{b-a+1}}
{(b-a+1)\lambda^{b-a+1}}}-c
\right|^p}{|x|^{\frac{pbN}{N-p}}}
e^{{-\frac{p|x|^{b-a+1}}
{(b-a+1)\lambda^{b-a+1}}}}
\mathrm{d}x
\nonumber\\&\quad
=
C(N,p,b)\lambda^{-(b-a+1)}
\inf_{c\in\mathbb{R}}
\int_{\mathbb{R}^N}
\frac{
\left|u-ce^{-\frac{|x|^{b-a+1}}
{(b-a+1)\lambda^{b-a+1}}}
\right|^p}{|x|^{pa}}
\mathrm{d}x
\nonumber\\&\quad
\ge C(N,p,b)\left(
\frac{\int_{\mathbb{R}^N}
\frac{|\nabla u|^p}{|x|^{pb}}\mathrm{d}x}
{\int_{\mathbb{R}^N}
\frac{|u|^p}{|x|^{pa}}\mathrm{d}x}
\right)^{\frac{1}{p}}
\inf_{c\in\mathbb{R},\lambda>0}
\int_{\mathbb{R}^N}
\frac{
\left|u-ce^{-\frac{|x|^{b-a+1}}
{(b-a+1)\lambda^{b-a+1}}}
\right|^p}{|x|^{pa}}
\mathrm{d}x.
\end{align*}
This completes the proof.
\end{proof}

\begin{rem}\label{rem-3.7}
Under $p\ge2$, the following assumptions also need hold to prove Theorem \ref{thm-2},
\begin{small}
\begin{equation}\label{4.16}
\begin{cases}
N\ge1,
\
p>1,
\
b-a+1>0,
\
b\le\frac{N-p}{p},\\
\sigma=\frac{p}{b-a+1}>0,\\
N-p>\varrho=pb\ge 0,\\
\theta=b-a+1\ge\frac{N-p-\varrho}{N-p},\\
pa=\frac{N\varrho}{N-p}=\frac{pbN}{N-p},
\end{cases}
\Longleftrightarrow
N>p\ge2,
\
0\le b<\frac{N-p}{p},
\
a=\frac{Nb}{N-p}.
\end{equation}
\end{small}
\end{rem}

\subsection{The stability results of the scale non-invariant $L^p$-CKN inequalities: proof of Theorems \ref{thm-3}, \ref{thm-4} and \ref{thm-5}}\label{sect-5}

\noindent In this section, we focus on the stability of the scale non-invariant $L^p$-CKN inequalities in three cases:
\begin{itemize}
[itemsep=0pt, topsep=0pt, parsep=0pt]

\item
$N>2$, $1<p<2$, $0\le b
<\frac{N-p}{p}$ and
$a=\frac{Nb}{N-p}$;

\item
$N>p\ge2$, $0\le b
<\frac{N-p}{p}$ and
$a=\frac{Nb}{N-p}$;

\item
$N>p\ge2$, $0\le b
<\frac{N-p}{p}$, $a<\frac{Nb}{N-p}$ and
$(p-1)a+b+1=\frac{pbN}{N-p}$.

\end{itemize}
Here notice that the above three situations correspond to Theorem \ref{thm-3}, Theorem \ref{thm-4} and Theorem \ref{thm-5},  respectively. We will start with proving Theorem \ref{thm-3}.

\begin{proof}
[\rm\textbf{Proof of Theorem \ref{thm-3}}]
From \eqref{4.1}, we get
\begin{align}\label{4.17}
&\int_{\mathbb{R}^N}
\frac{|\nabla u|^p}{|x|^{pb}}\mathrm{d}x +(p-1)\int_{\mathbb{R}^N} \frac{|u|^p}{|x|^{pa}}\mathrm{d}x -\left[N-(p-1)a-b-1\right] \int_{\mathbb{R}^N} \frac{|u|^p}{|x|^{(p-1)a+b+1}} \mathrm{d}x
\nonumber\\&\quad=
\int_{\mathbb{R}^N}
\frac{1}{|x|^{pb}}
\mathcal{G}_p\left(-u|x|^{b-a-1}x,
\nabla u\right)
\mathrm{d}x
\nonumber\\&\quad=
\int_{\Theta_{u,1}}
\frac{1}{|x|^{pb}}
\mathcal{G}_p\left(-u|x|^{b-a-1}x,
\nabla u\right)
\mathrm{d}x
+\int_{\Theta_{u,2}}
\frac{1}{|x|^{pb}}
\mathcal{G}_p\left(-u|x|^{b-a-1}x,
\nabla u\right)
\mathrm{d}x
\nonumber\\&\quad=I_1+I_2,
\end{align}
where $\Theta_{u,1}\cup\Theta_{u,2}
=\mathbb{R}^N$, and
\begin{align*}
\Theta_{u,1}
&:=\left\{x\in\mathbb{R}^N:
\left|u|x|^{b-a-1}x\right|
\le\left|\nabla u+u|x|^{b-a-1}x\right|
\right\},\\
\Theta_{u,2}&
:=\left\{x\in\mathbb{R}^N:
\left|u|x|^{b-a-1}x\right|
>\left|\nabla u+u|x|^{b-a-1}x\right|
\right\}.
\end{align*}
More precisely, $\Theta_{u,1}$ and $\Theta_{u,2}$ are Borel sets because of $u\in\mathcal{C}_0^\infty
(\mathbb{R}^N\setminus
\{\mathbf{0}\})\setminus\{0\}$, and it remains to apply Lemmas \ref{lem-2.5} and \ref{lem-2.7} to estimate $I_1$ and $I_2$. The proof of estimating $I_1$ and $I_2$ is similar to that of $\Pi_1$ and $\Pi_2$ (which is presented in the proof of Theorem \ref{thm-1}),
\begin{equation}\label{4.18}
I_1\ge
C(N,p,b)\inf_{c\in\mathbb{R}}
\int_{\Theta_{u,1}}
\frac{
\left|u-ce^{-\frac{|x|^{b-a+1}}
{b-a+1}}
\right|^p}{|x|^{\frac{pbN}{N-p}}}
\mathrm{d}x,
\end{equation}
and
\begin{equation}\label{4.19}
I_2\ge
C(N,p,b)
\inf_{c\in\mathbb{R}}
\int_{\Theta_{u,2}}
\frac{
\left||u|^{\frac{p-2}{2}}u
-|c|^{\frac{p-2}{2}}c
e^{-\frac{\frac{p}{2}|x|^{b-a+1}}
{b-a+1}}
\right|^2}{|x|^{\frac{N(pa+2b-2a)}{N-2}}}
\mathrm{d}x.
\end{equation}
Then, substituting \eqref{4.18} and \eqref{4.19} into \eqref{4.17}, it follows from \eqref{4.6} that
\begin{align*}
&\int_{\mathbb{R}^N} \frac{|\nabla u|^p}{|x|^{pb}}\mathrm{d}x +(p-1)\int_{\mathbb{R}^N} \frac{|u|^p}{|x|^{pa}}\mathrm{d}x -\left[N-(p-1)a-b-1\right] \int_{\mathbb{R}^N} \frac{|u|^p}{|x|^{(p-1)a+b+1}} \mathrm{d}x
\\&\quad
\ge
C(N,p,b)
\inf_{c\in\mathbb{R}}
\int_{\Theta_{u,1}}
\frac{
\left|u-ce^{-\frac{|x|^{b-a+1}}
{b-a+1}}
\right|^p}{|x|^{\frac{pbN}{N-p}}}
\mathrm{d}x
\\&\qquad
+C(N,p,b)
\inf_{c\in\mathbb{R}}
\int_{\Theta_{u,2}}
\frac{
\left||u|^{\frac{p-2}{2}}u
-|c|^{\frac{p-2}{2}}c
e^{-\frac{\frac{p}{2}|x|^{b-a+1}}
{b-a+1}}
\right|^2}{|x|^{\frac{N(pa+2b-2a)}{N-2}}}
\mathrm{d}x
\\&\quad
=
C(N,p,b)
\inf_{c\in\mathbb{R}}
\int_{\Theta_{u,1}}
\frac{
\left|u-ce^{-\frac{|x|^{b-a+1}}
{b-a+1}}
\right|^p}{|x|^{pa}}
\mathrm{d}x
\\&\qquad
+C(N,p,b)
\inf_{c\in\mathbb{R}}
\int_{\Theta_{u,2}}
\frac{
\left||u|^{\frac{p-2}{2}}u
-|c|^{\frac{p-2}{2}}c
e^{-\frac{\frac{p}{2}|x|^{b-a+1}}
{b-a+1}}
\right|^2}{|x|^{pa}}
\mathrm{d}x
\\&\quad\ge
C(N,p,b)
\inf_{c\in\mathbb{R}}
\int_{\mathbb{R}^N}
\frac{
\left||u|^{\frac{p-2}{2}}u
-|c|^{\frac{p-2}{2}}c
e^{-\frac{\frac{p}{2}|x|^{b-a+1}}
{b-a+1}}
\right|^2}{|x|^{pa}}
\mathrm{d}x.
\end{align*}
In view of the above arguments, the proof is completed.
\end{proof}

\begin{rem}
All the assumptions of Theorem \ref{thm-3} can be inferred from Remark \ref{rem-1} similarly.
\end{rem}

Now, we are ready to prove Theorem \ref{thm-4}.

\begin{proof}
[\rm\textbf{Proof of Theorem \ref{thm-4}}]
For all $q>0$ and $m,n\in\mathbb{R}$, there is a constant $C_q>0$ satisfying
\begin{equation}\label{4.20}
|m+n|^q\le C_q\left(|m|^q+|n|^q\right).
\end{equation}
Let $u:=\nu e^{-\frac{|x|^{b-a+1}}{b-a+1}}$, it follows from \eqref{2.2}, \eqref{4.1}, \eqref{4.20} and Lemma \ref{lem-2.2} (with $\varrho=pb$,  $\sigma=\frac{p}{b-a+1}$ and $\theta=b-a+1$) that
\begin{align*}
&\int_{\mathbb{R}^N}
\frac{\left|\nabla \left(u
-ce^{-\frac{|x|^{b-a+1}}{b-a+1}}\right)\right|^p}
{|x|^{pb}}\mathrm{d}x
+\int_{\mathbb{R}^N}
\frac{\left|u
-ce^{-\frac{|x|^{b-a+1}}
{b-a+1}}\right|^p}
{|x|^{pa}}\mathrm{d}x
\nonumber\\&\quad
=\int_{\mathbb{R}^N}
\frac{\left|\nabla \nu\cdot e^{-\frac{|x|^{b-a+1}}{b-a+1}}
-(\nu-c)e^{-\frac{|x|^{b-a+1}}
{b-a+1}}|x|^{b-a-1}x
\right|^p}
{|x|^{pb}}\mathrm{d}x
+\int_{\mathbb{R}^N}
\frac{|\nu-c|^p e^{-\frac{p|x|^{b-a+1}}{b-a+1}}}
{|x|^{pa}}\mathrm{d}x
\nonumber\\&\quad\le
C\int_{\mathbb{R}^N}
\frac{\left|\nabla \nu
\right|^p e^{-\frac{p|x|^{b-a+1}}{b-a+1}}}
{|x|^{pb}}\mathrm{d}x
+C\int_{\mathbb{R}^N}
\frac{|\nu-c|^p e^{-\frac{p|x|^{b-a+1}}{b-a+1}}}
{|x|^{pa}}\mathrm{d}x
\nonumber\\&\quad\le C\int_{\mathbb{R}^N}
\frac{\left|\nabla\nu
\right|^p e^{-\frac{p|x|^{b-a+1}}{b-a+1}}}
{|x|^{pb}}\mathrm{d}x
=C\int_{\mathbb{R}^N}
\frac{\left|\nabla
\left(ue^{\frac{|x|^{b-a+1}}
{b-a+1}}\right)\right|^p
e^{-\frac{p|x|^{b-a+1}}{b-a+1}}}
{|x|^{pb}}\mathrm{d}x
\nonumber\\&\quad=
C\int_{\mathbb{R}^N}
\frac{\left|\nabla u+u|x|^{b-a-1}x\right|^p}
{|x|^{pb}}\mathrm{d}x
\le C\int_{\mathbb{R}^N}
\frac{\mathcal{G}_p
\left(-u|x|^{b-a-1}x,\nabla u\right)}
{|x|^{pb}}\mathrm{d}x
\nonumber\\&\quad
=C\left\{\int_{\mathbb{R}^N} \frac{|\nabla u|^p}{|x|^{pb}}\mathrm{d}x +(p-1)\int_{\mathbb{R}^N} \frac{|u|^p}{|x|^{pa}}\mathrm{d}x -\left[N-(p-1)a-b-1\right] \int_{\mathbb{R}^N} \frac{|u|^p}{|x|^{(p-1)a+b+1}} \mathrm{d}x\right\},
\end{align*}
for some $C=C(N,p,b)>0$, which ends the proof.
\end{proof}

\begin{rem}
Using the same arguments as those of \eqref{4.16}, we know the source of the assumptions of Theorem \ref{thm-4}.
\end{rem}

We now conclude the discussion by presenting the proof of our final result (namely, Theorem \ref{thm-5}).

\begin{proof}
[\rm\textbf{Proof of Theorem \ref{thm-5}}]
By applying  \eqref{2.2}, \eqref{4.1} and Lemma \ref{lem-2.2} (with $\varrho=pb$,  $\sigma=\frac{p}{b-a+1}$ and $\theta=b-a+1$), we conclude that
\begin{align*}
\int_{\mathbb{R}^N}
\frac{\left|u
-ce^{-\frac{|x|^{b-a+1}}
{b-a+1}}\right|^p}
{|x|^{(p-1)a+b+1}}\mathrm{d}x
&=\int_{\mathbb{R}^N}
\frac{\left|ue^{\frac{|x|^{b-a+1}}
{b-a+1}}-c\right|^p
e^{-\frac{p|x|^{b-a+1}}{b-a+1}}}
{|x|^{(p-1)a+b+1}}\mathrm{d}x
\\&\le C(N,p,b)\int_{\mathbb{R}^N}
\frac{\left|\nabla
\left(ue^{\frac{|x|^{b-a+1}}
{b-a+1}}\right)\right|^p
e^{-\frac{p|x|^{b-a+1}}{b-a+1}}}
{|x|^{pb}}\mathrm{d}x
\\&=
C(N,p,b)\int_{\mathbb{R}^N}
\frac{\left|\nabla u+u|x|^{b-a-1}x\right|^p}
{|x|^{pb}}\mathrm{d}x
\\&\le C(N,p,b)\int_{\mathbb{R}^N}
\frac{\mathcal{G}_p
\left(-u|x|^{b-a-1}x,\nabla u\right)}
{|x|^{pb}}\mathrm{d}x
\\&=C(N,p,b)\bigg\{\int_{\mathbb{R}^N} \frac{|\nabla u|^p}{|x|^{pb}}\mathrm{d}x +(p-1)\int_{\mathbb{R}^N} \frac{|u|^p}{|x|^{pa}}\mathrm{d}x
\\&\quad
-\left[N-(p-1)a-b-1\right] \int_{\mathbb{R}^N} \frac{|u|^p}{|x|^{(p-1)a+b+1}} \mathrm{d}x\bigg\}.
\end{align*}
The proof is completed.
\end{proof}

\begin{rem}
Under the assumption $p\ge2$, we analyze the assumptions of Theorem \ref{thm-5}.
\begin{enumerate}[itemsep=0pt, topsep=1pt, parsep=0pt]

\item[(1)] We first analyze the source of the assumptions of Theorem \ref{thm-5}. In the proof process of Theorem \ref{thm-5}, in order to apply Lemma \ref{lem-4.1} and Lemma \ref{lem-2.2} (with $\varrho=pb$, $\sigma=\frac{p}{b-a+1}$ and $\theta=b-a+1$), there must hold:
\begin{equation*}
\begin{cases}
N\ge1,
\
p>1,
\
b-a+1>0,
\
b\le\frac{N-p}{p},\\
\sigma=\frac{p}{b-a+1}>0,\\
N-p>\varrho=pb\ge 0,\\
\theta=b-a+1\ge\frac{N-p-\varrho}{N-p}
=1-\frac{pb}{N-p},\\
(p-1)a+b+1=\frac{N\varrho}{N-p}
=\frac{pbN}{N-p},
\end{cases}
\end{equation*}
that is,
\begin{equation*}
N>p>1,
\
0\le b
<\frac{N-p}{p},
\
a\le\frac{Nb}{N-p},
\
(p-1)a+b+1=\frac{pbN}{N-p}.
\end{equation*}
Then, taking the intersection of the above set and $p\ge2$,
\[
N>p\ge2,
\
0\le b
<\frac{N-p}{p},
\
a\le\frac{Nb}{N-p},
\
(p-1)a+b+1=\frac{pbN}{N-p}.
\]
Moreover, if $a=\frac{Nb}{N-p}$, it derives from $(p-1)a+b+1=\frac{pbN}{N-p}$ that $(p-1)a+b+1=pa$, i.e., $b-a+1=0$. From Theorem \ref{thm-a}-(c),  there are no nontrivial functions such that $C(N,a,b)$ of \eqref{2-ckn-ineq} is achieved. Hence,
\begin{equation}\label{4.21}
N>p\ge2,
\
0\le b
<\frac{N-p}{p},
\
a<\frac{Nb}{N-p},
\
(p-1)a+b+1=\frac{pbN}{N-p}.
\end{equation}

\item[(2)] $b=\frac{N-p}{p+1}$ and $a=\frac{N}{p+1}-\frac{1}{p^2-1}$ for $N>p\ge2$ satisfy all the assumptions of \eqref{4.21}, which implies that the set presented in \eqref{4.21} is not empty.

\end{enumerate}
\end{rem}

\section*{Declarations}

\subsection*{Funding}
\noindent This paper was supported by the National Natural Science Foundation of China (No. 12371120).

\subsection*{Data availability statement}
\noindent No data was used for the research described in the article.

\subsection*{Conflict of interest}
\noindent The authors declare no conflict of interest.

\appendix

\section{Appendix: a technical inequality}\label{sectpls}

\noindent In this section, we focus on providing the detailed proof of \eqref{4.6}, which plays an important role in the proof of Theorem \ref{thm-1}.

\begin{lem}\label{lem-a}
For $m,n\in\mathbb{R}$ and $1<p<2$, the following inequality holds:
\begin{equation}\label{a.1}
\left||m|^{\frac{p-2}{2}}m
-|n|^{\frac{p-2}{2}}n\right|^2
\le 4|m-n|^p.
\end{equation}
In detail, if $\frac{m}{n}>0$,
\begin{equation}\label{a.2}
\left||m|^{\frac{p-2}{2}}m
-|n|^{\frac{p-2}{2}}n\right|^2
\le |m-n|^p;
\end{equation}
if $\frac{m}{n}<0$,
\begin{equation}\label{a.3}
\left||m|^{\frac{p-2}{2}}m
-|n|^{\frac{p-2}{2}}n\right|^2
\le 4|m-n|^p.
\end{equation}
\end{lem}

\begin{proof}[\rm\textbf{Proof}]
When $n=0$ or $m=0$, \eqref{a.1} obviously holds. Next, we consider $n\neq0$ and $m\neq0$. We divide the proof into the following two cases: $\frac{m}{n}>0$ and $\frac{m}{n}<0$.

\textbf{Case I: $\boldsymbol{\frac{m}{n}>0}$.} In this case, $\frac{m}{n}=\frac{|m|}{|n|}>0$ and
\begin{align*}
\eqref{a.2}\Longleftrightarrow
\left|
|m|^{\frac{p}{2}}
-|n|^{\frac{p}{2}}\right|^2
\le\Big||m|-|n|\Big|^p
&\Longleftrightarrow
\left|
\left(\frac{|m|}{|n|}
\right)^{\frac{p}{2}}
-1\right|^2
\le\left|\frac{|m|}{|n|}-1\right|^p
\\&\Longleftrightarrow
\left|
\left(\frac{m}{n}
\right)^{\frac{p}{2}}
-1\right|
\le\left|
\frac{m}{n}-1\right|^{\frac{p}{2}}
\\&\stackrel
{t=\frac{m}{n}}{\Longleftrightarrow}
\big|t^{\frac{p}{2}}-1\big|
\le\left|t-1\right|^{\frac{p}{2}},
\ t>0.
\end{align*}
Hence, it remains to show $\big|t^{\frac{p}{2}}-1\big|
\le\left|t-1\right|^{\frac{p}{2}}$ for all $t>0$.

$\bullet$ If $t=1$, $\big|t^{\frac{p}{2}}-1\big|
\le\left|t-1\right|^{\frac{p}{2}}$ holds obviously for $1<p<2$.

$\bullet$ If $t>1$, we see that $\big|t^{\frac{p}{2}}-1\big|
-\left|t-1\right|^{\frac{p}{2}}
=t^{\frac{p}{2}}-1
-\left(t-1\right)^{\frac{p}{2}}$. Let $f(t):
=t^{\frac{p}{2}}-1
-\left(t-1\right)^{\frac{p}{2}}$ for $1<p<2$. Then
\[
f'(t)=\frac{p}{2}
\left[t^{\frac{p-2}{2}}
-
\left(t-1\right)^{\frac{p-2}{2}}
\right]
=\frac{p}{2}
\left(t-1\right)^{\frac{p-2}{2}}
\left[
\left(\frac{t}{t-1}\right)^{\frac{p-2}{2}}
-1\right]<0,
\]
it shows that $f(t)<f(1)=0$ for all $t>1$. Thereby, for all $t>1$, $\big|t^{\frac{p}{2}}-1\big|
<\left|t-1\right|^{\frac{p}{2}}$.

$\bullet$ If $0<t\le\frac{1}{2}$, we get $\big|t^{\frac{p}{2}}-1\big|
-\left|t-1\right|^{\frac{p}{2}}
=1-t^{\frac{p}{2}}
-\left(1-t\right)^{\frac{p}{2}}$. Let $g(t):=1-t^{\frac{p}{2}}
-\left(1-t\right)^{\frac{p}{2}}$ for $1<p<2$. Notice that
\[
g'(t)=-\frac{p}{2}
\left[t^{\frac{p-2}{2}}
-\left(1-t\right)^{\frac{p-2}{2}}
\right]
=-\frac{p}{2}
\left(1-t\right)^{\frac{p-2}{2}}
\left[
\left(\frac{t}{1-t}
\right)^{\frac{p-2}{2}}
-1\right]
\ge0,
\]
then $g(t)\le g(\frac{1}{2})
=1-2^{\frac{2-p}{2}}<0$ for all $0<t\le\frac{1}{2}$. Therefore, for all $0<t\le\frac{1}{2}$, $\big|t^{\frac{p}{2}}-1\big|
<\left|t-1\right|^{\frac{p}{2}}$.

$\bullet$ If $\frac{1}{2}<t<1$, we get $g'(t)<0$, then $g(t)<g(\frac{1}{2})<0$ for all $\frac{1}{2}<t<1$. Thus, when $\frac{1}{2}<t<1$, $\big|t^{\frac{p}{2}}-1\big|
<\left|t-1\right|^{\frac{p}{2}}$.

In conclusion, $\big|t^{\frac{p}{2}}-1\big|
\le \left|t-1\right|^{\frac{p}{2}}$ for all $t>0$.

\textbf{Case II: $\boldsymbol{\frac{m}{n}<0}$.} In this case, $\frac{|m|}{|n|}>0$ and
\begin{align*}
\eqref{a.3}\Longleftrightarrow
\left|
|m|^{\frac{p}{2}}
+|n|^{\frac{p}{2}}\right|^2
\le4\Big||m|+|n|\Big|^p
&\Longleftrightarrow
\left|
\left(\frac{|m|}{|n|}
\right)^{\frac{p}{2}}
+1\right|^2
\le4\left|\frac{|m|}{|n|}
+1\right|^p
\\&\Longleftrightarrow
\left|\left(\frac{|m|}{|n|}
\right)^{\frac{p}{2}}
+1\right|
\le2\left|\frac{|m|}{|n|}
+1\right|^{\frac{p}{2}}
\\&\stackrel{t
=\frac{|m|}{|n|}}{\Longleftrightarrow}
\big|t^{\frac{p}{2}}+1\big|
\le2\left|t+1\right|^{\frac{p}{2}}
\\&\Longleftrightarrow
t^{\frac{p}{2}}+1
\le2\left(t+1\right)^{\frac{p}{2}},
\ t>0.
\end{align*}
In view of this, it suffices to prove $t^{\frac{p}{2}}+1
\le2\left(t+1\right)^{\frac{p}{2}}$ for all $t>0$.

Let $h(t):=t^{\frac{p}{2}}+1
-2\left(t+1\right)^{\frac{p}{2}}$ for $t>0$ and $1<p<2$. Observe that
\begin{align*}
h(t)&=\left[t^{\frac{p}{2}}
-
\left(t+1\right)^{\frac{p}{2}}\right]
+\left[1
-
\left(t+1\right)^{\frac{p}{2}}\right]
=\left(t+1\right)^{\frac{p}{2}}
\left\{
\left[
\left(
\frac{t}{t+1}\right)^{\frac{p}{2}}
-1\right]
+\left[
\left(\frac{1}{t+1}
\right)^{\frac{p}{2}}
-1\right]\right\}
<0.
\end{align*}
That is to say, $t^{\frac{p}{2}}+1
<2\left(t+1\right)^{\frac{p}{2}}$ for all $t>0$.
\end{proof}

\end{document}